\documentclass[a4paper,12pt]{amsart}
\usepackage{pdflscape}
\usepackage{amssymb}
\usepackage{amsmath}
\usepackage{array}
\usepackage{tikz}
\usetikzlibrary{arrows.meta}
\newtheorem{theorem}{Theorem}[section] 
\newtheorem{lemma}[theorem]{Lemma} 
\newtheorem{corollary}[theorem]{Corollary} 
\newtheorem{proposition}[theorem]{Proposition}
 
\theoremstyle{definition}
\newtheorem{definition}[theorem]{Definition} 
\newtheorem{remark}[theorem]{Remark}

\newcommand{\de}{\partial}
\newcommand{\slc}{\mathfrak{gl}_3}
\newcommand{\g}{\mathfrak{g}}
\newcommand{\F}{\mathbb F}

\newcommand{\C}{\mathbb F}
\newcommand{\Z}{\mathbb Z}
\title[Singular vectors in $E(3,8)$]{Degenerate Verma modules over the Lie superalgebra $E(3,8)$: a complete classification of singular vectors and morphisms}

\DeclareMathOperator{\diver}{div}
\DeclareMathOperator{\tr}{tr}
 
  \DeclareMathOperator{\Hom}{Hom}
\newcommand{\inlinewedge}{\textrm{\raisebox{0.6mm}{\footnotesize $\bigwedge$}}}

\author{Nicoletta Cantarini}\author{Fabrizio Caselli}\author{Victor Kac}
\address{Nicoletta Cantarini and Fabrizio Caselli, Dipartimento di matematica, Universit\`a di Bologna, Piazza di Porta San Donato 5, 40126 Bologna, Italy; Victor Kac, Math.\ Department, MIT, Cambridge, MA02139}

\email{nicoletta.cantarini@unibo.it}
\email{fabrizio.caselli@unibo.it}
\email{kac@math.mit.edu}

\begin{document}

\maketitle

\section{Introduction}
The list of infinite-dimensional linearly compact Lie superalgebras over an algebraically closed field of characteristic zero consists of 10 series and 5 exceptional Lie superalgebras: $E(1,6)$, $E(3,6)$, $E(3,8)$, $E(4,4)$, and $E(5,10)$ \cite{K}.
An explicit construction of the exceptional ones was given in \cite{CK2}, but, unfortunately, the construction of $E(3,8)$ given there contained a serious gap. This gap was filled in \cite{CCaK}, but the construction given there contained a typo, which was corrected in \cite{CaK}.

Using the construction of \cite{CK2}, Kac and Rudakov gave a classification of all degenerate (i.e., reducible) finite Verma modules over $E(3,8)$ announced in \cite[Theorem 4.5]{KR}, which turns out to be correct (see Theorem \ref{true} of the present paper). The reason for this is probably reflected in the fact that the missing elements in the construction of \cite{CK2} have very high degree. The morphisms, found in \cite{KR}, are depicted in Figure \ref{StartingFigure}.

In the present paper we prove the classification of reducible finite Verma modules over $E(3,8)$ and we find a complete list of morphisms between them. This is based on the classification of all singular vectors and the corresponding morphisms, given in Theorem \ref{CompleteMorphisms} (degrees 1,2,3,4), Theorems \ref{ClassDeg6} and \ref{ClassDeg6Morphisms} (degree 6), and Theorems \ref{ClassDeg7} and \ref{list} (degree 7). A complete proof of Theorem \ref{CompleteMorphisms}, announced in \cite{KR}, is given in Sections \ref{S4}--\ref{S7}. The absence of singular vectors of all other degrees is proved in Corollary \ref{bound} and a computer verification.

Our proofs are substantially simplified by the duality on the category of continuous $E(3,8)$-modules with discrete topology, established in \cite{CCK}.

All morphisms between degenerate finite Verma modules over $E(3,8)$ are depicted in Figure \ref{FinalFigure}, which provides an infinite set of complexes (Corollary \ref{V.}), that we call exceptional de Rham complexes for $E(3,8)$.

Throughout the paper the base field $\F$ is algebraically closed of characteristic zero.

\bigskip

\noindent{\bf Disclosure of AI use}. The mathematical results in this paper were developed without interaction with any AI agent. The TikZ code used to produce Figure~\ref{FinalFigure}  was generated with the assistance of ChatGPT (OpenAI) and subsequently checked by the authors. 
%That agent was also used for spell-checking.
Computations in finite Verma modules have been done with the assistance of a computer program developed by the authors in SAGE, without interactions with any AI agent. The code of this computer program is available upon request to the authors. 
\section{The setting}

We recall the structure of the Lie superalgebra $E(3,8)$ \cite{CK2, CCK, CaK}. 
%The description given in \cite{CCK} corrects the one in \cite{CK2}; a further correction was subsequently given in \cite{CaK}. 
The even and odd parts of $E(3,8)$ are given by
 \[E(3,8)_{\bar{0}}=W(3)+ (\Omega^0(3)\otimes \mathfrak{sl}_2)+d\Omega^1(3)\]
and 
\[E(3,8)_{\bar{1}}=(\Omega^0(3)^{-1/2}\otimes \C^2)+(\Omega^2(3)^{-1/2}\otimes\C^2),\]
where all sums are understood as direct sums of vector spaces, and
\begin{itemize}
    \item $W(3)$ is the Lie algebra of formal vector fields
     in the indeterminates $x_1, x_2, x_3$:
\[
W(3)=\Big\{\sum_{i=1}^3 P_i\frac{\partial}{\partial x_i}
\mid P_i\in\C[[x_1,x_2,x_3]]\Big\};\]
    \item $\mathfrak{sl}_2=\langle E=z_+\frac{\de}{\de z_-}, F=z_-\frac{\de}{\de z_+}, H=z_+\frac{\de}{\de z_+}-z_-\frac{\de}{\de z_-}\rangle$;
    \item $\C^2=\langle z_+,z_-\rangle$;
    \item $\Omega^k(3)$ is the space of formal differential $k$-forms in the variables $x_1,x_2,x_3$;
    \item $\Omega^k(3)^{\lambda}$, for $\lambda \in \C$, is the space $\Omega^k(3)$ endowed with a shifted action of $W(3)$: for every $X \in W(3)$ and $\omega\in \Omega^k(3)^\lambda$, the action of $X$ on $\omega$ is given by standard Lie derivative shifted by $\lambda \diver(X) \omega$.
\end{itemize} 

We let $X(\omega)$ denote the Lie derivative of a $k$-form $\omega\in \Omega^k(3)$ with respect to the vector field $X\in W(3)$.

We first describe the brackets between even elements. For $X,Y\in W(3)$, $f,g\in\Omega^0(3)$, $A,B\in \mathfrak{sl}_2$, $\omega_1, \omega_2\in d\Omega^1(3)$, the brackets are given by:

\begin{align*}
[X,Y]&=XY-YX-\frac{1}{2}d\big(\diver(X)\big)\wedge d\big(\diver(Y)\big),\\
[X,f\otimes A]&=X(f)\otimes A\\
[X,\omega_1]&=X(\omega_1)\\
[f\otimes A, g\otimes B]&=fg\otimes [A,B]+\tr(AB) df\wedge dg,\\
[f\otimes A,\omega_1]&=0,\\
[\omega_1,\omega_2]&=0.
\end{align*}
We next describe the brackets between even and odd elements. For $X\in W(3)$, $f\in\Omega^0(3)^{-1/2}$ $g\in\Omega^0(3)$, $v\in\C^2$, $A\in \mathfrak{sl}_2$, $\omega\in d\Omega^1(3)$,
$\sigma\in\Omega^2(3)^{-1/2}$, the brackets are given by
\begin{align*}
[X,f\otimes v]&=\big(X(f)-\frac{1}{2}\diver(X)f\big)\otimes v+\frac{1}{2}\big(d(\diver(X))\wedge df\big)\otimes v,\\
[X,\sigma \otimes v]&=(X(\sigma)-\frac{1}{2}\diver(X) \sigma) \otimes v\\
[g\otimes A,f\otimes v]&=(gf+dg\wedge df)\otimes A(v),\\
[g\otimes A, \sigma\otimes v]&=g\sigma\otimes A(v),\\
[\omega,f\otimes v]&=f\omega\otimes v,\\
[\omega, \sigma\otimes v]&=0.
\end{align*}
 
In order to describe the brackets between odd elements, it is convenient to introduce some notation. For every $2$-form $\omega\in\Omega^2(3)$, let $X_\omega\in W(3)$ be the unique vector field such that
\[
\iota_{X_\omega}(dx_1\wedge dx_2\wedge dx_3)=\omega.
\]
Similarly, for every $3$-form $\omega\in\Omega^3(3)$, let $f_\omega\in\Omega^0(3)$ be the unique function such that
\[
f_\omega\,dx_1\wedge dx_2\wedge dx_3=\omega.
\]

%Finally, we identify $W(3)$ with $\Omega^2(3)^{-1}$ and $\Omega^0(3)$ with $\Omega^3(3)^{-1}$. Besides, we identify $\Omega^2(3)^{-\frac{1}{2}}$ with $W(3)^{\frac{1}{2}}$ and we denote by $X_{\omega}$ the vector field corresponding to the $2$-form $\omega$ under this identification. We observe that if  $\omega_1, \omega_2\in \Omega^2(3)^{-1/2}$, then the Lie derivative $X_{\omega_1}(\omega_2)$ of $\omega_2$ along $X_{\omega_1}$ is a non-twisted 2-form. 

With this notation, the brackets between odd elements are given by the following formulas. For $\omega_1, \omega_2\in \Omega^2(3)^{-1/2}$, $f_1, f_2\in \Omega^0(3)^{-1/2}$, $u_1, u_2\in\C^2$, we have:
\begin{align*}
[f_1\otimes u_1, f_2\otimes u_2]&=(u_1\wedge u_2)(df_1\wedge df_2),\\
[f_1\otimes u_1, \omega_1\otimes u_2]&=(u_1\wedge u_2)(X_{f_1\omega_1}+df_1\wedge d(\diver(X_{\omega_1}))) -
\frac{1}{2}f_{f_1d\omega_1-\omega_1df_1}\otimes u_1\cdot u_2\\
[\omega_1\otimes u_1, \omega_2\otimes u_2]&=(u_1\wedge u_2)({X_{\omega_1}}(\omega_2)-\diver(X_{\omega_2})\omega_1).
\end{align*}
Here $u_1\cdot u_2$ denotes the symmetric product of $u_1$ and $u_2$, which we identify with an element of $\mathfrak{sl}_2$, while $u_1\wedge u_2$ denotes their skew-symmetric product, identified with a complex number. More precisely, we use the identifications
\[
E=\frac{1}{2} z_+^2,\qquad
F=-\frac{1}{2} z_-^2,\qquad
H=-z_+\cdot z_-,
\qquad
z_+\wedge z_-=1.
\]

We consider the $\mathbb{Z}$-grading on $E(3,8)$  defined by 
$$\deg x_i=-\deg \de_{i}=2,\,\,\deg E=\deg F=\deg H=0, \,\, \deg dx_i=2, \,\, \deg z_+=\deg z_-=-3,$$
for $i=1,2,3$.

We shall simplify the notation
by setting  for $i,j=1,2,3$ and $\epsilon \in \{+,-\}$:
\[
\begin{aligned}
&\partial_i=\frac{\partial}{\partial x_i},\qquad
\partial_\epsilon=\frac{\partial}{\partial z_\epsilon},\
&d_{ij}=dx_i\wedge dx_j,\qquad
d_i^\epsilon=x_i\otimes z_\epsilon,\qquad
d^\epsilon=1\otimes z_\epsilon.
\end{aligned}
\]

We shall also find it convenient to fix a total ordering of the odd elements of negative degree. Accordingly, we set
\[
n_1=d_1^+,\quad n_2=d_2^+,\quad n_3=d_3^+,\quad
n_4=d_1^-,\quad n_5=d_2^-,\quad n_6=d_3^-,
\quad n_7=d^+,\quad n_8=d^-.
\]

The above $\Z$-grading on $E(3,8)$ is an irreducible, consistent $\mathbb{Z}$-grading of depth 3, known as the principal grading of $E(3,8)$. Its degree-zero component is isomorphic to $\mathfrak{sl}_3\oplus \mathfrak{sl}_2\oplus\C$ and is spanned by
the elements $x_i\de_{j}$, $E$, $F$ and $H$. The element $Y=2\sum_{i=1}^3{x_i\de_i}$ turns out to be a grading element of $E(3,8)$. 

We denote by $V(a,b,c,t)$ the finite-dimensional irreducible $\g_0$-module with highest weight $(a,b,c,t)$ with respect to the Cartan elements
\[h_{12}=x_1\de_1-x_2\de_2, \quad h_{23}=x_2\de_2-x_3\de_3, \quad H, \quad Y.\]
Recall that $a,b,c\in \Z_{\geq 0}$, $t\in \F$.
Besides,
\begin{align*}E(3,8)_{-1}&=\langle d_i^+, d_i^- ~|~ i=1,2,3\rangle=\langle n_i ~|~ i=1,\ldots,6\rangle,\\
E(3,8)_{-2}&=\langle \de_{i}~|~ i=1,2,3\rangle,\\
E(3,8)_{-3}&=\langle d^+, d^-\rangle=\langle n_7,n_8\rangle.
\end{align*}
The component $E(3,8)_{-1}$ generates the negative part of the grading and, as an $E(3,8)_0$-module, is isomorphic to \[\C^3\boxtimes \C^2\boxtimes \C.\]
Similarly, $E(3,8)_{1}$ generates the positive part of the grading and is given by
\[E(3,8)_{1}=\langle x_ix_j\otimes z_\epsilon, d_{ij}\otimes z_\epsilon\,|\, i,j=1,2,3, \epsilon=+,-\rangle.\]
As an $E(3,8)_0$-module, it decomposes into two irreducible components, namely,  $E(3,8)_{1}\cong (S^2\C^3\boxtimes \C^2\boxtimes \C)\oplus (\inlinewedge ^2\C^3\boxtimes \C^2\boxtimes \C)$. The lowest weight vectors in the two components are, up to scalar factors, the elements $x_3^2\otimes z_-$ and $d_{23}\otimes z_-$, respectively.

\medskip

Let us denote by $\g$ the Lie superalgebra $E(3,8)$ endowed with its principal grading:
 \[\g=\prod_{j\geq -3}\g_j.\] 
 
Given a finite-dimensional irreducible $\g_0$-module $V$, we regard it as
a $\g_{\geq 0}$-module by letting $\g_+=\prod_{j>0}\g_j$ act trivially, and define
$$M(V)=U(\g)\otimes_{U(\g_{\geq 0})}V.$$
The module $M(V)$ is called the finite Verma module associated with $V$. It carries  a natural $\g$-module structure induced by the left multiplication. Moreover, as $\g_0$-module, we have \[M(V)\cong U(\g_{<0})\otimes_{\C}V, \qquad\g_{<0}=\g_{-1} \oplus \g_{-2}\oplus \g_{-3}.\]
A finite Verma module is said to be non-degenerate if it is irreducible and degenerate otherwise. Let us set $U_-=U(\g_{<0})$. We change the sign of the grading in $U_-$, so that it becomes a $\Z_{\geq 0}$-graded space.
\begin{definition}
We say that a non-zero element $w\in M(V)$ is  {\em homogeneous} of degree $d$ if $w\in (U_-)_{d}\otimes V$.
\end{definition}
Degenerate Verma modules can be described in terms of morphisms between Verma modules.
A morphism $\varphi: M(V)\rightarrow M(W)$ can be associated with an element 
$\Phi\in U_{-}\otimes \Hom(V,W)$ as follows. Writing
\[\Phi=\sum_iu_i\otimes \theta_i, \quad u_i\in U_-, \quad \theta_i\in \Hom(V,W),
\]
we define
$$\varphi(u\otimes v)=u\Phi(v),$$
for $u\in U_-$ and $v\in V$, where \[\Phi(v)=\sum_iu_i\otimes \theta_i(v).\]
We  say that $\varphi$ (or $\Phi$) is a morphism of degree $d$ if $u_i\in (U_-)_d$ for every $i$.

\medskip

The following proposition characterizes morphisms between finite Verma modules.

\begin{proposition}\cite{R}\label{morphisms}
Let $\varphi: M(V)\rightarrow M(W)$ be the linear map associated with the element $\Phi\in U_{-}\otimes \Hom(V,W)$.
Then $\varphi$ is a morphism of $\g$-modules if and only if the following conditions
hold:
\begin{itemize}
\item[(a)] $\g_0.\Phi=0$;
\item[(b)] $X\varphi(v)=0$ for every $X\in \g_1$ and for every $v\in V$.
\end{itemize} 
\end{proposition}

Moreover, if $M(V)$ is a finite Verma module and condition (a) holds, it is enough to verify  condition (b) for a set of generators of $\g_1$ as a $\g_0$-module and for a highest weight vector in $V$.

Kac and Rudakov in \cite[Theorem 4.5]{KR} announced  a classification of degenerate Verma modules for $E(3,8)$. They construct  certain morphisms between finite Verma modules of degree 1, 2, 3, and 4, thereby proving that all modules in their list are indeed degenerate. They also announced that the proof of the completeness of the list would be published elsewhere, but such a proof has never appeared.
The main goal of this paper is to prove that the list of degenerate finite Verma modules given by Kac and Rudakov is indeed complete. To this end, we provide a classification of all morphisms between finite Verma modules. Our classification strictly extends the one given by Kac and Rudakov: we discover two new families of morphisms of degree 6 and two new morphisms of degree 7. 

We can now state the classification of degenerate Verma modules.
\begin{theorem}\label{true} Let $V$ be a finite dimensional irreducible $\g_0$-module. Then $M(V)$ is degenerate if and only if $V$ is isomorphic to one of the following:
\begin{itemize}
    \item $V_A(a,c)=V(a,0,c,-4a+3c)$, for $a,c\in\Z_{\geq 0}$;
    \item $V_B(a,c)=V(a,0,c,-4a-3c-6)$, for $a,c\in\Z_{\geq 0}$;
    \item $V_C(b,c)=V(0,b,c,4b+3c+12)$, for $b,c\in\Z_{\geq 0}$;
    \item $V_D(b,c)=V(0,b,c,4b-3c+6)$, for $b,c\in\Z_{\geq 0}$ and $(b,c)\neq (0,0)$;
\end{itemize}
\end{theorem}
To describe the morphisms between Verma modules it is convenient to realize the modules $V_A(a,c), V_B(a,c), V_C(b,c), V_D(b,c)$ as follows. For all $a,b,c\geq 0$ we let 
\begin{itemize}
   \item $V_A(a,c) =S^a(x_1,x_2,x_3)\boxtimes S^c(z_+,z_-)$, with $Y$ acting by $-4a+3c$;
   \item $V_B(a,c)=S^a(x_1,x_2,x_3)\boxtimes S^c(\de_+,\de_-)$, with $Y$ acting by $-4a-3c-6$;
   \item  $V_C(b,c)=S^b(\de_1,\de_2,\de_3)\boxtimes S^c(z_+,z_-)$, with $Y$ acting by $4b+3c+12$;
   \item $V_D(b,c)=S^b(\de_1,\de_2,\de_3)\boxtimes S^c(\de_+,\de_-)$, with $Y$ acting by $4b-3c+6$.
\end{itemize}
In particular, for $i=1,2,3$, we have that the action of $\de_i$ on $V_A$ and $V_B$ is given by differentiation, whereas on $V_C$ and $V_D$ is given by multiplication. Similarly, for $\epsilon=+,-$,  the action of $\de_\epsilon$ on $V_A$ and $V_C$ is given by differentiation, whereas on $V_B$ and $V_D$ it is given by multiplication.

For notational convenience we introduce the following notation: 
\begin{align*}
D_0&=d_1^-d_2^-d_3^-,\\
D_1&=d_1^+d_2^-d_3^-+d_1^-d_2^+d_3^-+d_1^-d_2^-d_3^+,\\ 
D_2&=d_1^+d_2^+d_3^-+d_1^+d_2^-d_3^++d_1^-d_2^+d_3^+, \\ 
D_3&=d_1^+d_2^+d_3^+.
\end{align*}

For $a,c\in\Z_{\geq 0}$, we set $M_A(a,c)=M(V_A(a,c))$, and we similarly define $M_B(a,c)$, $M_C(b,c)$, and $M_D(b,c)$. The following theorem summarizes the morphisms between finite Verma modules obtained by Kac and Rudakov.
\begin{theorem}\label{CompleteMorphisms} The following is a list of morphisms between Verma modules (see Figure \ref{StartingFigure} for an illustration).
\begin{enumerate}
    \item In degree 1 the element \[\sum_{i=1}^3\sum_{\epsilon=+,-}d_i^\epsilon \otimes \de_i\de_{\epsilon}\] defines the morphisms 
    \begin{align*}
        \varphi[1_A]&:M_A(a,c)\rightarrow M_A(a-1,c-1), \,\,\textrm{for }a,c\geq 1;\\
        \varphi[1_B]&:M_B(a,c)\rightarrow M_B(a-1,c+1), \,\,\textrm{for }a\geq 1, c\geq 0;\\
        \varphi[1_C]&:M_C(b,c)\rightarrow M_C(b+1,c-1), \,\,\textrm{for }b\geq 0,c\geq 1;\\
        \varphi[1_D]&:M_D(b,c)\rightarrow M_D(b+1,c+1),\,\, \textrm{for }b,c\geq 0.
    \end{align*}
    \item In degree 2 the element
    \[
    \sum_{i,j=1}^{3}d_i^+d_j^-\otimes \de_i \de_j
    \]
    defines the morphisms
    \begin{align*}
        \varphi[2_{AB}]&:M_A(a,0)\rightarrow M_B(a-2,0),\,\,\textrm{for }a\geq 2;\\
        \varphi[2_{CD}]&: M_C(0,b)\rightarrow M_D(0,b+2),\,\,\textrm{for }b\geq 0.\\
    \end{align*}
    \item In degree 3 the element
    \[
    \sum_{k=0}^3 D_k \otimes \de_+^k \de_-^{3-k}
    \]
    defines the morphisms
    \begin{align*}
       \varphi[3_{AC}]&:M_A(0,c)\rightarrow M_C(0,c-3),\,\,\textrm{for }c\geq 3;\\
        \varphi[3_{BD}]&: M_B(0,c)\rightarrow M_D(0,c+3),\,\,\textrm{for }c\geq 0.\\ 
    \end{align*}
    \item In degree 4 the element
    \[
    \sum_{k=0}^2 \sum_{i=1}^3 d_i^- D_{k+1}\otimes \de_i \de_+^{k}\de_-^{2-k}
    \]
    %\color{magenta}
    %We can also write
    %\[
    %\sum_{k=0}^2 \sum_{i=1}^3 d_i^+ D_{k}\otimes \de_i \de_+^{k}\de_-^{2-k}
    %\]
    %\color{black}
    defines the morphism
    \[
    \varphi[4_{AD}]:M_A(0,2)\rightarrow M_D(1,0)
    \]
    and the element
    \[
    \sum_{k=0}^2 \sum_{i=1}^3  D_{k+1} d_i^-\otimes \de_i \de_+^{k}\de_-^{2-k}
    \]
    defines the morphism
    \[
    \varphi'[4_{AD}]:M_A(1,0)\rightarrow M_D(0,2).
    \]
\end{enumerate}
\end{theorem}
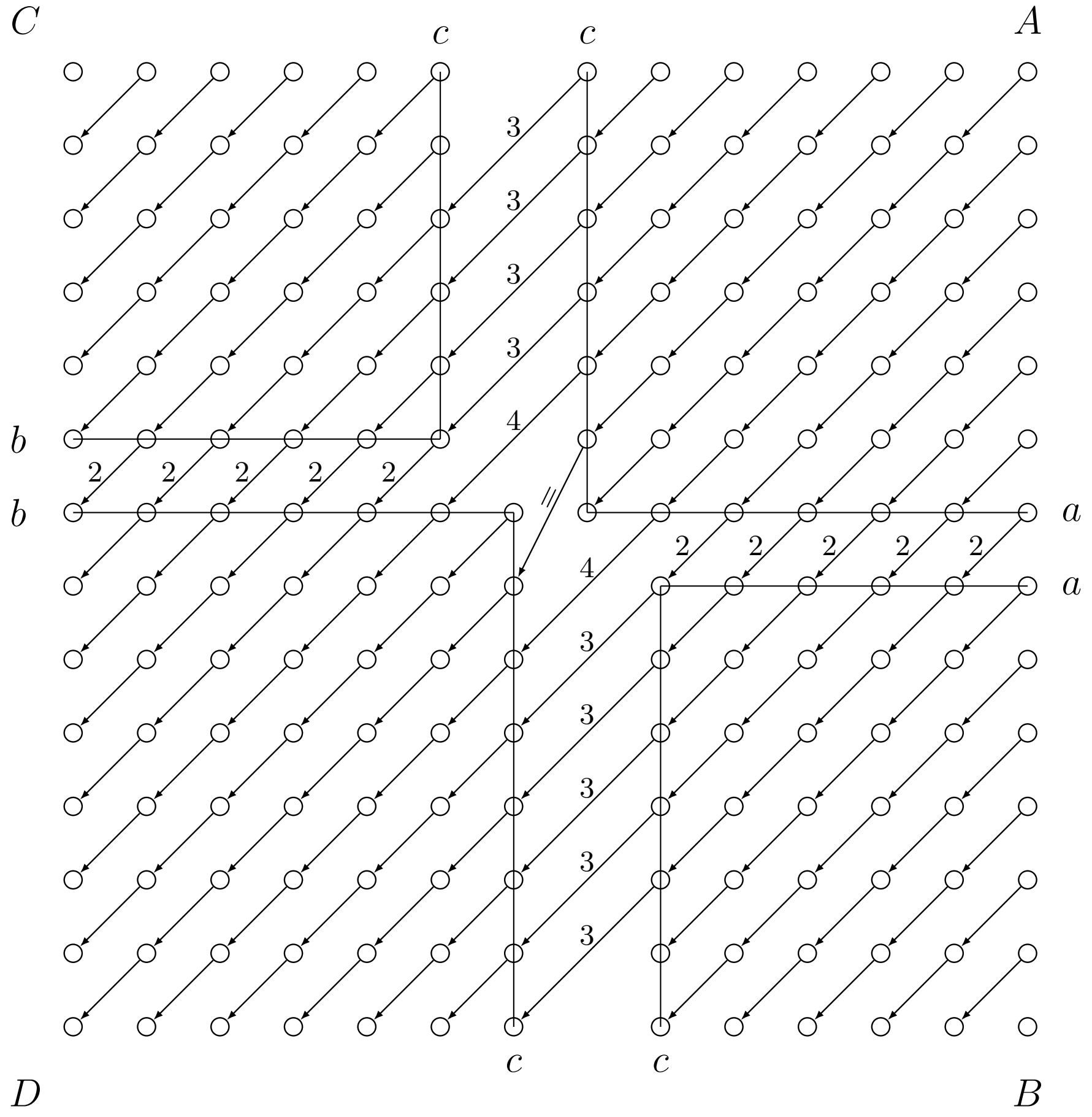
\begin{figure}
\begin{tikzpicture}[scale=1.1,
	x=0.72cm,
	y=0.72cm,
	>=Latex,
	dot/.style={
		circle,
		draw,
		fill=white,
		minimum size=5.5pt,
		inner sep=0pt,
		line width=.45pt
	},
	arr/.style={
		-{Latex[length=3pt,width=2.2pt]},
		line width=.45pt
	},
	darr/.style={
		-{Latex[length=3pt,width=2.2pt]},
		dash pattern=on 4pt off 4pt,
		line width=.55pt
	},
	barr/.style={
		-{Latex[length=4pt,width=3pt]},
		line width=1pt
	}
	]
	
	% =========================================================
	% QUADRANTE C
	% =========================================================
	
	\foreach \i in {0,...,5}{
		\foreach \j in {0,...,5}{
			\node[dot] (C\i\j) at (\i,\j) {};
		}
	}
	
	% frecce interne: da NE verso SW
	\foreach \i in {1,...,5}{
		\pgfmathtruncatemacro{\im}{\i-1}
		\foreach \j in {1,...,5}{
			\pgfmathtruncatemacro{\jm}{\j-1}
			\draw[arr] (C\i\j) -- (C\im\jm);
		}
	}
	
	\draw (C00.center) -- (C50.center);
	\draw (C50.center) -- (C55.center);

	% =========================================================
	% QUADRANTE A
	% =========================================================
	
	\foreach \i in {0,...,6}{
		\foreach \j in {0,...,6}{
			\node[dot] (A\i\j) at ({7+\i},\j-1) {};
		}
	}
	
	% frecce interne: da NE verso SW
	\foreach \i in {1,...,6}{
		\pgfmathtruncatemacro{\im}{\i-1}
		\foreach \j in {1,...,6}{
			\pgfmathtruncatemacro{\jm}{\j-1}
			\draw[arr] (A\i\j) -- (A\im\jm);
		}
	}
	
	\draw (A00.center) -- (A60.center);
	\draw (A00.center) -- (A06.center);

	% =========================================================
	% QUADRANTE D
	% =========================================================
	
	\foreach \i in {-1,...,5}{
		\foreach \j in {0,...,7}{
			\node[dot] (D\i\j) at (\i+1,-1-\j) {};
		}
	}
	
	% frecce interne
	\foreach \i in {0,...,5}{
		\pgfmathtruncatemacro{\im}{\i-1}
		\foreach \j in {0,...,6}{
			\pgfmathtruncatemacro{\jm}{\j+1}
			\draw[arr] (D\i\j) -- (D\im\jm);
		}
	}
	
	\draw (D50.center) -- (D57.center);
	\draw (D-10.center) -- (D50.center);

	% =========================================================
	% QUADRANTE B
	% =========================================================
	
	\foreach \i in {0,...,5}{
		\foreach \j in {0,...,6}{
			\node[dot] (B\i\j) at ({8+\i},-2-\j) {};
		}
	}
	
	% frecce interne
	\foreach \i in {1,...,5}{
		\pgfmathtruncatemacro{\im}{\i-1}
		\foreach \j in {0,...,5}{
			\pgfmathtruncatemacro{\jm}{\j+1}
			\draw[arr] (B\i\j) -- (B\im\jm);
		}
	}
	
	\draw (B00.center) -- (B50.center);
	\draw (B00.center) -- (B06.center);

	% =========================================================
	% FRECCE TRA QUADRANTI
	% =========================================================
	
	% A -> C
	
	%\draw[arr] (A02) -- (C50);
	\draw[arr] (A03) -- (C50);
	%\draw[arr] (A03) -- (C51);
	\draw[arr] (A04) -- (C51);
	%\draw[arr] (A04) -- (C52);
	\draw[arr] (A05) -- (C52);
	%\draw[arr] (A05) -- (C53);
	\draw[arr] (A06) -- (C53);
	%\draw[arr] (A06) -- (C54);
	% B -> D
	%\draw[arr] (B00) -- (D52);
	\draw[arr] (B00) -- (D53);
	%\draw[arr] (B01) -- (D53);
	\draw[arr] (B01) -- (D54);
	%\draw[arr] (B02) -- (D54);
	\draw[arr] (B02) -- (D55);
	%\draw[arr] (B03) -- (D55);
	\draw[arr] (B03) -- (D56);
	%\draw[arr] (B04) -- (D56);
	\draw[arr] (B04) -- (D57);
	%\draw[arr] (B05) -- (D57);
	
	%C-->D
	\draw[arr] (C50) -- (D30);
	\draw[arr] (C40) -- (D20);
	\draw[arr] (C30) -- (D10);
	\draw[arr] (C20) -- (D00);
	\draw[arr] (C10) -- (D-10);
	
	%A-->B
	\draw[arr] (A20) -- (B00);
	\draw[arr] (A30) -- (B10);
	\draw[arr] (A40) -- (B20);
	\draw[arr] (A50) -- (B30);
	\draw[arr] (A60) -- (B40);

	% A -> D
        %\draw[arr] (A10) -- (D51);
       % \draw[arr] (A01) -- (D40);
	\draw[arr] (A02) -- (D40);
	\draw[arr] (A10) -- (D52);
	\draw[arr] (A01) -- (D51);
	
	% =========================================================
	% ETICHETTE
	% =========================================================
	
	\node at (-0.65,5.7) {$C$};
	\node at (13.0,5.7) {$A$};
	
	\node at (-0.65,-8.9) {$D$};
	\node at (13.0,-8.9) {$B$};
	
	\node at (5,5.5) {$c$};
	\node at (7,5.5) {$c$};
	
	\node at (-0.75,0) {$b$};
	\node at (-0.75,-1) {$b$};
	
	\node at (13.6,-1) {$a$};
	\node at (13.6,-2) {$a$};
	
	\node at (6,-8.5) {$c$};
	\node at (8,-8.5) {$c$};
	
	% etichette gradi
	\foreach \i in {1,...,5}{
		\node at (\i-0.7,-0.45) {\scriptsize{2}};
	}

	\foreach \i in {1,...,5}{
		\node at (\i+7.3,-1.45) {\scriptsize{2}};
	}
	
	\foreach \j in {1,...,4}{
		\node at (6,\j+0.25) {\scriptsize{3}};
	}
	
	%\foreach \j in {1,...,5}{
	%	\node at (6.4,\j-0.1) {\scriptsize{6}};
	%}

     %   \node at (6.4, -0.1) {\scriptsize{7}};
    %    \node at (6.6, -1.5) {\scriptsize{7}};
    
	\foreach \j in {1,...,5}{
		\node at (7,-\j-1.75) {\scriptsize{3}};
	}
	
	%\foreach \j in {1,...,6}{
	%	\node at (7.4,-\j-1.1) {\scriptsize{6}};
	%}
	
	\node at (6,0.25) {\scriptsize{4}};
	\node at (7,-1.75) {\scriptsize{4}};
	\node[rotate=64] at (6.5,-0.8){\scriptsize{=}};

	\node at (1.5,-11) {$M_A(a,c)=M(a,0,c,-4a+3c)$;};
	\node at (10,-11) {$M_B(a,c)=M(a,0,c,-4a-3c-6)$};
	\node at (1.7,-12) {$M_C(b,c)=M(0,b,c,4b+3c+12);$};
	\node at (9.9,-12) {$M_D(b,c)=M(0,b,c,4b-3c+6)$.};
	
\end{tikzpicture}
\caption{\label{StartingFigure} Morphisms between finite Verma modules obtained  by Kac and Rudakov}
\end{figure}

Note that $(U_-)_{d}$ is finite-dimensional. We denote its dimension by $m_d$. For example, $m_1=6$, $m_2=18$, $m_3=40$, $m_4=78$. 
In the following sections, we shall use  the following basis of $U_-$: 
\[
\mathcal B=
\left\{
\de_1^{a_1}\de_2^{a_2}\de_3^{a_3}
n_{i_1}n_{i_2}\cdots n_{i_k}
:\,
a_1,a_2,a_3\in\Z_{\geq 0},\quad
1\leq i_1<i_2<\cdots<i_k\leq 8
\right\}.
\]

%where $n_i=d_i^+$ if $i=1,2,3$, $n_i=d_{i-3}^-$ if $i=4,5,6$, $n_7=1\otimes e_1$ and $n_8=1\otimes e_2$,
We denote by $\mathcal B_d$ the subset of $\mathcal B$ consisting of the elements of degree $d$; thus, $\mathcal B_d$
is a basis of $(U_-)_{d}$. Every vector $m$ in $M(V)_d$ can be uniquely expressed as 
\[
m=\sum_{u\in \mathcal B_d}u\otimes v_u,\]
where $v_u\in V$. We call $v_u\in V$ the \emph{coefficient} of $u$ in $m$. In particular, $m=0$ if and only if all its coefficients vanish.

\begin{definition}\label{singular} A vector $w\in M(V)$ is called  {\em singular} if it satisfies the following conditions:
\begin{itemize}
\item[(i)] $x_i\partial_{i+1}w=0$ for every $i=1,2$; \label{fields}
\item[(ii)] $Ew=0$;\label{E}
\item[(iii)] $Xw=0$ for every $X\in \g_1$.
\end{itemize}
\end{definition}
We observe that every homogeneous component of positive degree of a singular vector is again singular. 
The same holds for its weight components. Henceforth, unless otherwise specified, we shall assume that singular vectors are  homogeneous both in degree and in weight.

To classify singular vectors of degree $d$, we consider a
vector
\[
w=\sum_{u\in\mathcal B_d}u\otimes v_u
\]
and apply to it all the elements $X$ of a linear basis of
$\mathfrak g_1$. We then expand each $Xw$ with respect to the
basis $\mathcal B_{d-1}$ and require all its coefficients to
vanish.

Each of these coefficients is a linear combination of vectors
$v_u$ and of vectors of the form $xv_u$, where
$x\in\mathfrak g_0$ and $u\in\mathcal B_d$.

By considering only those equations whose coefficients involve
the vectors $v_u$ and no terms of the form $xv_u$, we obtain
a homogeneous linear system, which we call the
\emph{basic homogeneous linear system} in degree $d$.
 For example, the basic linear homogeneous system in degree 1 is empty, whereas in degree 2 it consists of the following equations: 
\[
\begin{cases}
	v_{\de_1}=0,\\
	v_{n_2n_3}=0,\\
	v_{n_5n_6}=0,\\
	v_{n_3n_5}-v_{n_2n_6}=0.
\end{cases}
\]

%{\color{blue}\begin{lemma}\label{singsing}Let $0\neq w\in M(V)$. Then $w$ is an $E(3,8)$-singular vector if and only if
%	\begin{itemize}
%		\item $(x_1\de_2)w=0=(x_2\de_3)w$;
%		\item $E w=0$; 
%		\item $(x_3^2\otimes e_2)w=0$;
%		\item $(d_{23}\otimes e_2)w=0$.
%	\end{itemize}
%end{lemma}
%\begin{proof}
%This is a consequence of Definition \ref{singular} and of the description of the principal grading of $\g$, since $x_3^2\otimes e_2$ and $d_{23}\otimes e_2$ are lowest weight vectors in the two irreducible components of $\g_1$.
%\end{proof}}

We recall that a finite Verma module $M(V)$ is degenerate if and only if it contains a singular vector of positive degree \cite[Proposition 3.3]{CC}.

 If $V$ is an irreducible finite dimensional $\g_{0}$-module of highest weight $\lambda=(a,b,c,t)$ we also denote the Verma module $M(V)$ by $M(\lambda)$ or $M(a,b,c,t)$. 

For a weight $\lambda=(a,b,c,t)$, we denote by $\lambda^{\vee}=(b,a,c,-t+6)$ its (shifted) dual.

We recall the following duality for finite Verma modules, established in \cite{CCK} in a general setting of continuous $E(3,8)$-modules with discrete
topology. 
%Recall that for a $\g_0$-module $F(\lambda)$ we denote by $\lambda^*$ the highest weight of $F(\lambda)^*$.
\begin{theorem}\label{duality}There exists a contravariant functor on the category of finite Verma modules $M(\lambda)\rightarrow M(\lambda^{\vee})$. More precisely, if $\varphi:M(\lambda)\rightarrow M(\mu)$ is a  $\g$-modules morphism of degree $d$, then there exists a dual morphism $\varphi^\vee:M(\mu^{\vee})\rightarrow M(\lambda^{\vee})$ of the same degree. Equivalently, if $M(\mu)$ contains a singular vector of degree $d$ and weight $\lambda$, then $M(\lambda^\vee)$ contains a singular vector of degree $d$ and weight $\mu^\vee$.
\end{theorem}

The classification of morphisms between finite Verma modules is equivalent to the classification of singular vectors in these modules \cite[Proposition 3.5]{CC}. Indeed, if $M(\mu)$ contains a singular vector $w$ of degree
$d$ and weight $\lambda$, and $s$ is a highest-weight vector
of degree $0$ and weight $\lambda$ in $M(\lambda)$, then
there exists a unique morphism
\[
\varphi:M(\lambda)\longrightarrow M(\mu)
\]
of degree $d$ such that $\varphi(s)=w$.

Conversely, every morphism
$\varphi:M(\lambda)\to M(\mu)$ is uniquely determined by
the image of a highest-weight vector of $M(\lambda)$, which is necessarily a singular vector of weight $\lambda$.
If $w$ is a singular vector which determines a morphism $M(\lambda)\rightarrow M(\mu)$ we denote by $w^\vee$ the singular vector (defined up to a scalar factor) which determines the dual morphism $M(\mu^\vee)\rightarrow M(\lambda^\vee)$.
\begin{definition}
	If $w$ is a singular vector of weight $\lambda$ in $M(\mu)$, we define $\operatorname{step}(w)=\mu-\lambda$.
\end{definition}
\begin{remark} \label{rem_dual} If $w$ and $w^\vee$ are dual singular vectors and \[\operatorname{step}(w)=(\alpha,\beta,\gamma, \tau),\] then \[\operatorname{step}(w^\vee)=(-\beta,-\alpha,-\gamma, \tau).\] 
\end{remark}
%\color{blue}
%\begin{remark}\label{dual}
%Let $\varphi: M(V)\rightarrow M(W)$ be a linear map of degree $d$ associated to an element $\Phi\in U_-\otimes \Hom(V,W)$ that satisfies condition (a) of Proposition \ref{morphisms}. Then there exists a $\g_0$-morphism $\psi: (U_-)_d^*\rightarrow \Hom(V,W)$ such that
%$\Phi= \sum_i u_i\otimes \psi(u_i^*)$ where $\{u_i, i\in I\}$ is any basis of $(U_-)_d$ and $\{u_i^*, i\in I\}$ is the corresponding dual basis.
%\end{remark}
%\color{black}

\section{A bound on the degree of singular vectors}
Let $\g=E(3,8)$.
The main goal of this section is to prove that every singular vector for $\g$ in a finite Verma module has degree at most 14. We consider the Lie subalgebra $\mathfrak h =W(3)+ d\Omega^1(3)$ of $\g$.

We denote by $\widetilde U_j$ the subspace of
$U_-$ spanned by all products of
$j$ elements of $\mathfrak g_{-1}\oplus\mathfrak g_{-3}$.

Let $V$ be a finite-dimensional $\mathfrak g_0$-module.
We consider the following filtration of $M(V)$ by vector
subspaces. For every $k\geq 0$, we set
\[
M_k(V)=
\F[\partial_1,\partial_2,\partial_3]\otimes
\sum_{j\leq k}\widetilde U_j\otimes V.
\]
By construction, we have
\[
M_0(V)\subset M_1(V)
\subset\cdots\subset M_8(V)=M(V),
\]
since $\dim(\mathfrak g_{-1}\oplus\mathfrak g_{-3})=8$.

The following proposition shows that this filtration is
compatible with the action of the even part of $\mathfrak g$.
\begin{proposition}\label{s5mod}
	For every $k=0,1,\ldots,8$ the subspace $M_k(V)$ is a $\g_{\bar 0}$-submodule of $M(V)$.
\end{proposition}
\begin{proof}
	It is enough to prove that, for every $X\in \g_{\bar 0}$, every $1\leq j\leq k$, every $u_1,\ldots,u_j\in \{\g_{-1}\oplus \g_{-3}\}$, and every $v\in V$, we have
	\begin{equation}\label{s5mod1}
		Xu_1\cdots u_j\otimes v \in M_j(V).
	\end{equation}
We shall prove simultaneously that
	\begin{equation}
		\label{s5mod2} [X,u_1]u_2\cdots u_j\otimes v\in M_j(V).
	\end{equation}
	The proof is by a double induction on $j$ and $\deg X$. 

    For $j=1$, the inclusion in \eqref{s5mod2} follows
directly from the grading, and hence so does
\eqref{s5mod1}.

If $\deg X=-2$, both inclusions follow directly from
the grading. We may therefore assume that $j\geq 2$
and $\deg X\geq 0$.

        We have
\[
Xu_1\cdots u_j\otimes v
=
[X,u_1]u_2\cdots u_j\otimes v
+
u_1Xu_2\cdots u_j\otimes v.
\]
The second summand belongs to $M_j(V)$ by induction
on $j$, and therefore \eqref{s5mod1} follows from
\eqref{s5mod2}. Moreover,
\[
\begin{aligned}
[X,u_1]u_2\cdots u_j\otimes v
={}&-u_2[X,u_1]u_3\cdots u_j\otimes v\\
&+[[X,u_1],u_2]u_3\cdots u_j\otimes v.
\end{aligned}
\]
The first summand belongs to $M_j(V)$ by induction
on $j$, while the second belongs to $M_j(V)$ by
induction on $\deg X$. This completes the proof.
	
\end{proof}
By Proposition~\ref{s5mod}, the filtration
\[
\{0\}=M_{-1}(V)\subseteq
\mathbb F[\partial_1,\partial_2,\partial_3]\otimes V
=M_0(V)\subseteq\cdots\subseteq M_8(V)=M(V)
\]
is a filtration by $\mathfrak h$-submodules.
We set
\[
N_k(V)=M_k(V)/M_{k-1}(V),
\qquad k=0,\ldots,8.
\]

Recall that we denote by $n_1,\ldots,n_8$ the following elements:
\[
d_1^+,d_2^+,d_3^+,d_1^-,d_2^-,d_3^-,d^+,d^-,
\]
respectively.

\begin{proposition}\label{NkV}
	For every $k=0,\ldots,8$, we have an isomorphism of $\C$-vector spaces
	\[
	N_k(V)\cong \F[\de_1,\de_2,\de_3]\otimes \bigoplus_{p_1< \cdots < p_k} \F n_{p_1}\cdots n_{p_k}\otimes V
	\]
	as $\C$-vector spaces. 
\end{proposition}
\begin{proof}
	For every $p_1,\ldots,p_k \in \{1,\ldots,8\}$ and every permutation $\sigma\in S_k$ of the indices $\{1,\ldots,k\}$, we have 
	\begin{equation}\label{permform}n_{p_1}n_{p_2}\cdots n_{p_k}-\varepsilon(\sigma) n_{p_{\sigma(1)}}\cdots n_{p_{\sigma(k)}}\in M_{k-1}(V),
	\end{equation}
    where $\varepsilon(\sigma)$ denotes the sign of $\sigma$.
	It follows that $N_k(V)$ is generated, as $\F[\de_1,\de_2,\de_3]$-module, by the elements \[n_{p_1}\cdots n_{p_k}\otimes v, \qquad p_1< \cdots < p_k, \quad v\in V.\] The claimed isomorphism now follows from the  Poincar\'e--Birkhoff--Witt theorem for $U_-$.
\end{proof}	

The subspace
\[
F_k(V)=\bigoplus_{p_1< \cdots < p_k}\F n_{p_1}\cdots n_{p_k}\otimes V
\]
of $N_k(V)$ has a particularly simple structure.
\begin{proposition}\label{NkVisVerma}
	The subspace $F_k(V)$ of $N_k(V)$ is a $\slc$-module annihilated by $\mathfrak h_{>0}$.
	Moreover, $N_k(V)$ is the finite Verma module for $\mathfrak h$ induced by $F_k(V)$, i.e.
	\[
	N_k(V)\cong M(F_k(V)).
	\]
\end{proposition}
\begin{proof}
The subspace $F_k(V)$ is a $\slc$-module
because $\mathfrak g_{-1}\oplus\mathfrak g_{-3}$ is a
$\mathfrak g_0$-module and, by \eqref{permform}, the
products of its elements are alternating in the
associated quotient.

Furthermore, $F_k(V)$ is annihilated by
$\mathfrak h_{>0}$, since
$[\mathfrak h_2,\mathfrak g_{-3}]=0$.
The second assertion follows from Proposition~\ref{NkV}
and the first part of the proof.
	
\end{proof}
\begin{lemma}\label{lemhgl3}
Let $V$ be a $\slc$-module, and consider the $\mathfrak h$-Verma module $M(V)$ induced by $V$. If \[w=\sum_{I}\partial_I\otimes v_I\in M(V)\] 
is a singular vector, then the same vector is singular in  the $W(3)$-Verma module induced by $V$.
\end{lemma}
\begin{proof}
This follows from the fact that, in $E(3,8)$, the bracket $[X,\partial_i]$ is simply the commutator $X\de_i-\de_i X$, since $\diver(\de_i)=0$.
\end{proof}
These results allow us to relate singular vectors for $E(3,8)$ to singular vectors for $W(3)$. We shall use the following classical classification theorem due to Rudakov \cite{R1}.
\begin{theorem}\label{thmw3}
Let $a,b\in\Z_{\geq 0}$ and $F_t(a,b)$ be the irreducible $\slc$-module with highest weight $(a,b)$ with respect to $h_{12}$, $h_{23}$ and central charge $t$ with respect to $\sum_ix_i\de_i$. Let $w$ be a singular vector of positive degree in $M(F_t(a,b))$. Then, up to a  scalar factor, $w$ is one of the following vectors:
\begin{itemize}
    \item $\de_3\otimes 1 \in M(F_0(0,0))$;
    \item $\de_2\otimes x_3^*-\de_3 \otimes x_2^* \in M(F_{-1}(0,1))$;
    \item $\de_1\otimes x_1+\de_2 \otimes x_2 +\de_3 \otimes x_3 \in M(F_{-2}(1,0)).$
\end{itemize}
\end{theorem}

We can now establish the first bound on the degree of singular vectors for $E(3,8)$.
\begin{corollary}\label{bound}
	Let $M(V)$ be a finite Verma module for $E(3,8)$ and $w\in M(V)$ be a singular vector. Then $w$ has degree at most 14.
\end{corollary}
\begin{proof}
	Let $k$ be the smallest integer such that $w\in M_k(V)$, and let $\bar w$ be the nonzero image of $w$ in $N_k(V)$. Two cases may occur.
	\begin{enumerate}
	\item If $\bar w$ is a highest-weight vector  in $F_k(V)$, then  $w$ has degree at most 12. This is because there are six generators of degree $1$ and two
generators of degree $3$ in
$\mathfrak g_{-1}\oplus\mathfrak g_{-3}$, hence every element
of $F_k(V)$ has degree at most
\[
6+2\cdot 3=12.
\]
	\item Otherwise, $\bar w$ is a singular vector of positive degree in the $\mathfrak h$-Verma module $N_k(V)$. In particular, by Lemma \ref{lemhgl3} and Theorem \ref{thmw3}, its degree as a singular vector
for $\mathfrak h$ is 2 (since in $E(3,8)$ we have $\deg \de_i=2$). Hence $w$ has degree at most $12+2=14$.
\end{enumerate}
\end{proof}

\section{Degree 1}\label{S4}
In this section, we classify all singular vectors of degree 1 and  show that they are precisely those obtained by Kac and Rudakov.  Let \[w=\sum_{i=1}^3(d_i^+\otimes v_i^+)+\sum_{i=1}^3(d_i^-\otimes v_i^-),\] 
 where $v^+_i, v^-_i \in V$, be a singular vector of degree 1 and weight $(a,b,c,t)$ in a finite Verma module $M(V)$.

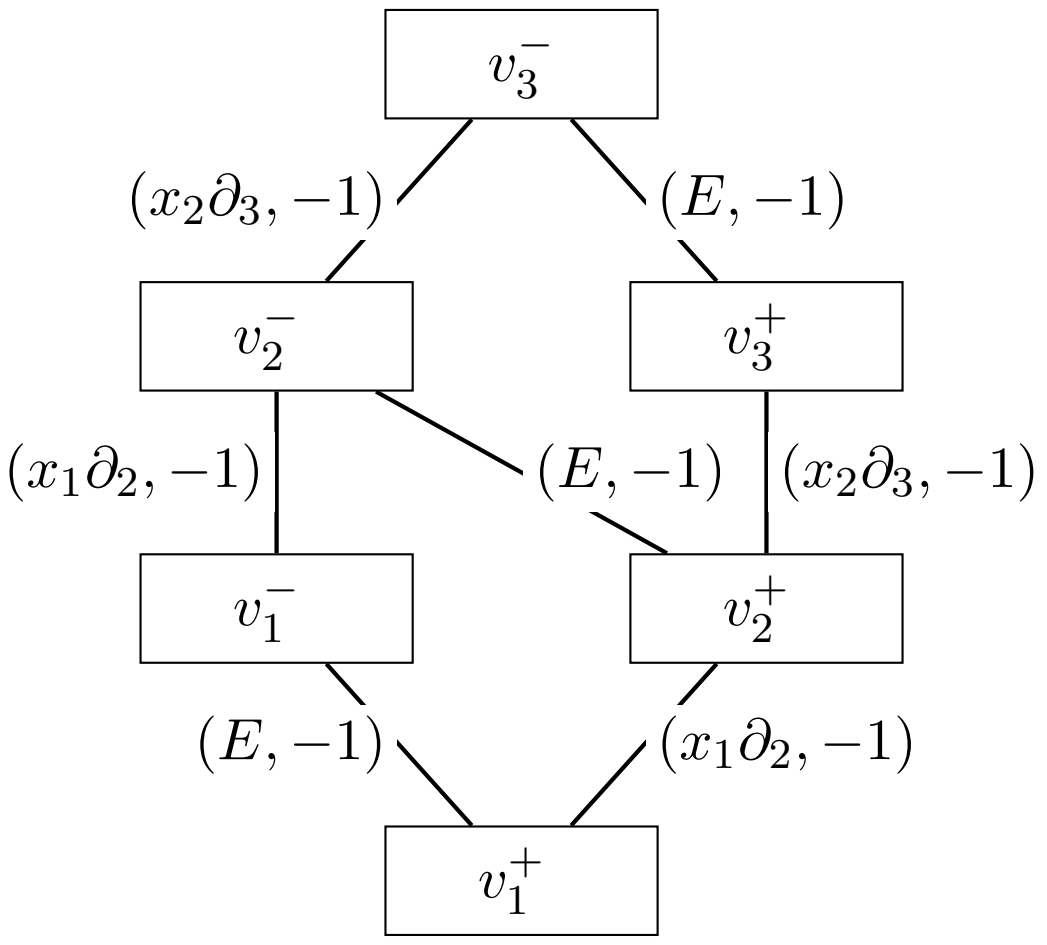
\begin{figure}[h]
\begin{tikzpicture}[scale=0.9,
	box/.style={
		draw,
		minimum width=1.8cm,
		minimum height=0.6cm,
		align=center,
		font=\small
	},
	edge/.style={
		draw,
		thick
	},
	lab/.style={
		font=\small,
		fill=white,
		inner sep=2pt
	}
	]
% ============================================================
% TERZO DIAGRAMMA
% ============================================================

\node[box] (C1) at (18,7)
{$v_3^-$};

\node[box] (C2) at (16.2,5)
{$v_2^-$
	%		$({a-1}\!+\!1,\!{b+1},\!{c+2}\!+\!2)$
};

\node[box] (C3) at (19.8,5)
{$v_3^+$
	%		$({a}\!+\!1,\!{b},\!{c})$
};

\node[box] (C4) at (16.2,3)
{$v_1^-$
	%		$({a},\!{b}\!+\!1,\!{c}\!+\!2)$
};

\node[box] (C5) at (19.8,3)
{$v_2^+$
	%		$({a}\!-\!1,\!{b}\!+\!1,\!{c})$
};

\node[box] (C6) at (18,1)
{$v_1^+$
	%		$({a},\!{b}\!-\!1,\!{c})$
};

% edges
\draw[edge] (C1) -- node[lab,left] {$(x_2\partial_3,-1)$} (C2);
\draw[edge] (C1) -- node[lab,right] {$(E,-1)$} (C3);

\draw[edge] (C2) -- node[lab,left] {$(x_1\partial_2,-1)$} (C4);
\draw[edge] (C2) -- node[lab,right] {$(E,-1)$} (C5);

\draw[edge] (C3) -- node[lab,right] {$(x_2\partial_3,-1)$} (C5);

\draw[edge] (C4) -- node[lab,left] {$(E,-1)$} (C6);
\draw[edge] (C5) -- node[lab,right] {$(x_1\partial_2,-1)$} (C6);

\end{tikzpicture}

\caption{\label{FigureDegree1} The action of $x_1\de_2$, $x_2\de_3$ and $E$ on the  degree-one components.} 

\end{figure}

The conditions $x_i\de_{i+1}w=0$, for $i=1,2$, and $Ew=0$ (see Definition \ref{singular}) are equivalent to the following system:

\begin{enumerate}
\item $\left\{
\begin{array}{l}
(x_1\de_2) v_1^{\pm}+ v_2^{\pm}=0,\\
(x_1\de_2) v_i^{\pm}=0\,\, i=2,3;\\
\end{array}\right.$
\item $\left\{
\begin{array}{l}
(x_2\de_3) v_2^{\pm}+ v_3^{\pm}=0,\\
(x_2\de_3) v_i^{\pm}=0\,\, i=1,3;\\
\end{array}\right.$

\item $\left\{
\begin{array}{l}
E v_i^{+}+ v_i^{-}=0,\\
E v_i^{-}=0,
\end{array}\right.$
for $i=1,2,3$. 
\end{enumerate}
The action of $x_1\de_2,x_2\de_3$ and $E$ on the vectors $v_i^{\pm}$ is summarized by the directed graph shown in Figure \ref{FigureDegree1}, where we assume that all edges are oriented from bottom to top. 

More precisely, an edge from
a vertex $v$ to a vertex $v'$, labelled by $(X,c)$,
indicates that $X(v)=cv'$, where
$X\in\{x_1\partial_2,x_2\partial_3,E\}$ and
$c\in\C$. If no edge corresponding to $X$ starts
at $v$, then $X(v)=0$.

\bigskip

\begin{table}
\begin{tabular}[]{|l|l|l|}
	\hline
Element	 & Weight & Step\\
	\hline
	\hline
	$v_1^+$ & $(a-1,b,c-1,t+1)$&$(-1,0,-1,1)$\\
	\hline 
	$v_2^+$ & $(a+1,b-1,c-1,t+1)$&$(1,-1,-1,1)$\\
	\hline 
	$v_3^+$ & $(a,b+1,c-1,t+1)$&$(0,1,-1,1)$\\
	\hline
	$v_1^-$ & $(a-1,b,c+1,t+1)$&$(-1,0,1,1)$\\
	\hline
	$v_2^-$ & $(a+1,b-1,c+1,t+1)$&$(1,-1,1,1)$\\
	\hline
	$v_3^-$ & $(a,b+1,c+1,t+1)$&$(0,1,1,1)$.\\
	\hline 
\end{tabular}
\bigskip
\caption{Weights in degree one.}
\label{weightsDegree1}
\end{table}

The weights of the vectors $v_i^{\pm}$, together with their corresponding steps, are listed in Table \ref{weightsDegree1}.
In particular, the table gives $\operatorname{step}(w)$
when one of the vectors $v_i^\pm$ is a highest-weight
vector of $V$.

Let $w^\vee \in M(V')$ be the singular vector dual to $w$. We can construct the analogous table by replacing
$(a,b,c,t)$ with the weight of $w^\vee$. To distinguish
the vectors associated with $w$ from those associated
with $w^\vee$, we write $v_i^\pm(w)$ and
$v_i^\pm(w^\vee)$, respectively.

By Remark~\ref{rem_dual}, if $v_1^+(w)$ (respectively,
$v_2^+(w)$ or $v_3^+(w)$) is a highest-weight vector
of $V$, then $v_3^-(w^\vee)$ (respectively,
$v_2^-(w^\vee)$ or $v_1^-(w^\vee)$) is a highest-weight
vector of $V'$, and conversely.

Therefore, up to duality, we may assume that the
highest-weight vector of $V$ is one of
$v_1^+$, $v_2^+$ or $v_1^-$. Under this assumption,
the vectors $v_2^-$, $v_3^+$ and $v_3^-$ must vanish,
for weight reasons.

\begin{proposition}\label{sing1simple} Let  
	\[
	w=d_1^+\otimes v_1^++d_2^+\otimes v_2^++d_1^-\otimes v_1^-\in M(V)
	\]
be a singular vector in a finite Verma module. Then, up to a scalar factor, $w$ is one of the following:
\begin{itemize}
	\item $w[1_A]=d_1^+\otimes x_1^{a-1} z_+^{c-1}\in M_A(a-1,c-1)$ for all $a,c\geq 1$; 
	\item $w[1_B]=d_1^+\otimes x_1^{a-1} \de_{+}\de_{-}^c+d_1^-\otimes x_1^{a-1}\de_{-}^{c+1}\in M_B(a-1,c+1)$ for all $a\geq 1$, $c\geq 0$.
\end{itemize}
\end{proposition}
\begin{proof}
	Applying $x_1x_3\otimes z_-$ to $w$, we obtain $(a+b+1)v_2^+=0$ and so $v_2^+=0$ since $a,b\geq 0$.
    
	Next, applying $x_2x_3\otimes z_+$, $x_2x_3\otimes z_-$, $d_{23}\otimes z_+$ and $d_{23}\otimes z_-$ to $w$ we obtain the following equations
	
	\[\begin{cases}
		bv_1^-=0,\\
		bv_1^+=0,\\
		(-4a - 2b - 3c -t - 6)v_1^-=0,\\
		(4a + 2b - 3c + t)v_1^+-6Fv_1^-=0.
	\end{cases}
	\]
Since $w\neq 0$ and $v_2^+=0$, at least one of
$v_1^+$ and $v_1^-$ is nonzero. Hence $b=0$.

If $v_1^-\neq 0$ then  $v_1^-$ is a highest-weight vector and $t=-4a-3c-6$. Moreover, the last equation gives 
\[v_1^+=-\frac{1}{c+1}Fv_1^-.\] 
This yields the singular vector $w[1_B]$.

If $v_1^-=0$, then $v_1^+\neq 0$, and we obtain $b=0$ and $t=-4a+3c$. This gives the singular vector $w[1_A]$.
\end{proof}

%\begin{enumerate}
%\item[(i)] $(x_1\de_2) w=0$;
%\item[(ii)] $(x_2\de_3) w=0$;
%\item[(iii)] $E w=0$;
%\item[(iv)] $(x_3^2\otimes e_2) w=0$;
%\item[(v)] $(d_{23}\otimes e_2) w=0$.
%\end{enumerate}

\begin{theorem}\label{classDegree1}
The following is a complete list, up to scalar factors, of singular vectors of degree 1 in finite Verma modules, together with the corresponding morphisms.
\begin{itemize}
\item $w[1_A]=d_1^+\otimes x_1^{a-1} z_+^{c-1}\in M_A(a-1,c-1)$ for all $a,c\geq 1$. It defines the morphism
$\varphi[1_A]: M_A(a,c)\rightarrow M_A(a-1,c-1)$;\vspace{2mm}
\item $w[1_B]=d_1^+\otimes x_1^{a-1}\de_{+}\de_{-}^c+d_1^-\otimes x_1^{a-1}\de_{-}^{c+1}\in M_B(a-1,c+1)$ for all $a\geq 1$, $c\geq 0$. It defines the morphism
$\varphi[1_B]: M_B(a,c)\rightarrow M_B(a-1,c+1)$;\vspace{2mm}
\item $w[1_C]=\sum_{i=1}^3 d_i^+\otimes \de_i \de_3^b z_+^{c-1}%=d_1^+\de_1\de_3^{b}+d_2^+\de_2\de_3^{b}+d_3^+\de_3^{b+1})\otimes z_1^{c-1}
\in M_C(b+1,c-1)$ for all $b\geq 0$, $c\geq 1$. It defines the morphism
$\varphi[1_C]: M_C(b,c)\rightarrow M_C(b+1,c-1)$;\vspace{2mm}
\item $w[1_D]=\sum_{i=1}^3 \sum_{\epsilon=+,-}d_i^\epsilon\otimes \de_i \de_3^b \de_\epsilon \de_-^c%(d_1^+\de_1\de_3^{b}+d_2^+\de_2\de_3^{b}+d_3^+\de_3^{b+1})\otimes \de_{z_1}\de_{z_2}^c+(d_1^-\de_1\de_3^{b}+d_2^-\de_2\de_3^{b}+d_3^-\de_3^{b+1})\otimes \de_{z_2}^{c+1}
\in M_D(b+1,c+1)$ for all $b, c\geq 0$. It defines the morphism
$\varphi[1_D]: M_D(b,c)\rightarrow M_D(b+1,c+1)$.
\end{itemize}
%For every $X=A,B,C,D$, we have $\varphi[1_X]=\sum_{i=1}^3\delta_i\de_i$, where
%\[\delta_i=d_i^+\otimes \de_{z_1}+d_i^-\otimes \de_{z_2}.\]
\end{theorem}

\begin{proof}
By Proposition~\ref{sing1simple} and the preceding
discussion, every singular vector of degree $1$ is,
up to a scalar factor, one of $w[1_A]$, $w[1_B]$
or their duals.

A direct computation shows that $w[1_C]$ and
$w[1_D]$ are singular vectors. By the duality
described in Remark~\ref{rem_dual}, they are
respectively the duals of $w[1_B]$ and $w[1_A]$.
This completes the proof.
\end{proof}

\section{Degree 2}
In this section, we classify singular vectors of degree $2$
and the corresponding morphisms between finite Verma modules.

Let $w_{\mathrm{sing}}(2)$ be a singular vector of degree $2$
and weight $(a,b,c,t)$:
\[
w_{\mathrm{sing}}(2)
=
\sum_{i=1}^3 \de_i\otimes v_i
+
\sum_{1\leq i<j\leq 6}
n_i n_j\otimes v_{\overline{ij}}
\in M(V),
\]
where, as before, $n_i=d_i^+$ for $i=1,2,3$ and
$n_i=d_{i-3}^-$ for $i=4,5,6$, and $V$ is a
finite-dimensional irreducible $\g_0$-module.

\begin{figure}[h]
		%\rotatebox{90}{%
			
\begin{tikzpicture}[scale=0.9,
	box/.style={
		draw,
		minimum width=1.3cm,
		minimum height=0.6cm,
		align=center,
		font=\small
	},
	edge/.style={
		draw,
		thick
	},
	lab/.style={
		font=\tiny,
		fill=white,
		inner sep=2pt
	}
	]
	
	% ============================================================
	% PRIMO DIAGRAMMA
	% ============================================================
	
	\node[box] (A1) at (0,7)
	{$v_{\overline{36}}$
		%$({a},{b}\!+\!2,{c})$
		};
	
	\node[box] (A2) at (0,5)
	{$v_{\overline{26}}+v_{\overline{35}}$
%		$({a}\!+\!1,{b},{c})$
};
	
	\node[box] (A3) at (-1.8,3)
	{$v_{\overline{25}}$
%		$({a}\!+\!2,{b}\!-\!2,{c})$
};
	
	\node[box] (A4) at (1.8,3)
	{$v_{\overline{16}}+v_{\overline{34}}$
%		$({a}\!-\!1,{b}\!+\!1,{c})$
};
	
	\node[box] (A5) at (0,1)
	{$v_{\overline{15}}+v_{\overline{24}}$
%		$({a},{b}\!-\!1,{c})$
};
	
	\node[box] (A6) at (0,-1)
	{$v_{\overline{14}}$
%		$({a}\!-\!2,{b},{c})$
};
	
	% edges
	\draw[edge] (A1) -- node[lab,right] {$(x_2\partial_3,-2)$} (A2);
	
	\draw[edge] (A2) -- node[lab,left] {$(x_2\partial_3,-1)$} (A3);
	\draw[edge] (A2) -- node[lab,right] {$(x_1\partial_2,-1)$} (A4);
	
	\draw[edge] (A3) -- node[lab,left] {$(x_1\partial_2,-2)$} (A5);
	\draw[edge] (A4) -- node[lab,right] {$(x_2\partial_3,-1)$} (A5);
	
	\draw[edge] (A5) -- node[lab,right] {$(x_1\partial_2,-1)$} (A6);

	% ============================================================
	% SECONDO DIAGRAMMA
	% ============================================================
	
	\node[box] (B1) at (7.5,7)
	{$v_{\overline{56}}$
%		$({a}\!+\!1,{b},{c}\!+\!2)$
};
	
	\node[box] (B2) at (7.5,5)
	{$v_{\overline{46}}$
%		$({a}\!-\!1,{b}\!+\!1,{c}\!+\!2)$
};
	
	\node[box] (B3) at (10.5,5)
	{$v_{\overline{26}}-v_{\overline{35}}$
%		$({a}\!+\!1,{b},{c})$
};
	
	\node[box] (B4) at (7.5,3)
	{$v_{\overline{45}}$
%		$({a},\!{b}\!-\!1,\!{c}\!+\!2)$
};
	
	\node[box] (B5) at (10.5,3)
	{$v_{\overline{16}}-v_{\overline{34}}$
%		$({a}\!-\!1,\!{b}\!+\!1,\!{c})$
};
	
	\node[box] (B6) at (13.5,3)
	{$v_{\overline{23}}$
%		$({a}\!+\!1,\!{b}\!+\!1,\!{c}\!-\!2)$
};
	
	\node[box] (B7) at (7.5,1)
	{$v_{\overline{15}}-v_{\overline{24}}$
	%	$({a},{b}\!-\!1,{c})$
	};
	
	\node[box] (B8) at (10.5,1)
	{$v_{\overline{13}}$
%		$({a}\!-\!1,{b}\!+\!1,{c}\!-\!2)$
};
	
	\node[box] (B9) at (9,-1)
	{$v_{\overline{12}}$
%		$({a},{b}\!-\!1,{c}\!-\!2)$
};
	
	% edges
	\draw[edge] (B1) -- node[lab] {$(x_1\partial_2,\!-1)$} (B2);
	\draw[edge] (B1) -- node[lab,right] {$(E,-2)$} (B3);
	
	\draw[edge] (B2) -- node[lab] {$(x_2\partial_3,-1)$} (B4);
	\draw[edge] (B2) -- node[lab] {$(E,-2)$} (B5);
	
	\draw[edge] (B3) -- node[lab] {$(x_1\partial_2,-1)$} (B5);
	\draw[edge] (B3) -- node[lab,right] {$(E,-1)$} (B6);
	
	\draw[edge] (B4) -- node[lab] {$(E,-2)$} (B7);
	\draw[edge] (B5) -- node[lab] {$(x_2\partial_3,-1)$} (B7);
	\draw[edge] (B5) -- node[lab] {$(E,-1)$} (B8);
	
	\draw[edge] (B6) -- node[lab,right] {$(x_1\partial_2,-1)$} (B8);
	
	\draw[edge] (B7) -- node[lab,left] {$(E,-1)$} (B9);
	\draw[edge] (B8) -- node[lab,right] {$(x_2\partial_3,-1)$} (B9);

	% ============================================================
	% TERZO DIAGRAMMA
	% ============================================================

	\node[box] (C3) at (4.5,5)
	{$v_{1}$
%		$({a}\!+\!1,\!{b},\!{c})$
};
	
	\node[box] (C5) at (4.5,3)
	{$v_{2}$
%		$({a}\!-\!1,\!{b}\!+\!1,\!{c})$
};
	
	\node[box] (C6) at (4.5,1)
	{$v_{3}$
%		$({a},\!{b}\!-\!1,\!{c})$
};
	
	% edges
	%\draw[edge] (C1) -- node[lab,left] {$(x_1\partial_2,-1)$} (C2);
	\draw[edge] (B1) -- node[lab,left] {$(E,+1)$} (C3);
	
	%\draw[edge] (C2) -- node[lab,left] {$(x_2\partial_3,-1)$} (C4);
	\draw[edge] (B2) -- node[lab,left] {$(E,-1)$} (C5);
	
	\draw[edge] (C3) -- node[lab,left] {$(x_1\partial_2,+1)$} (C5);
	
	\draw[edge] (B4) -- node[lab,left] {$(E,+1)$} (C6);
	\draw[edge] (C5) -- node[lab,left] {$(x_2\partial_3,+1)$} (C6);
	
\end{tikzpicture}

\caption{\label{FigureDegree2} The action of $x_1\de_2$, $x_2\de_3$ and $E$ in degree 2.} 

\end{figure}
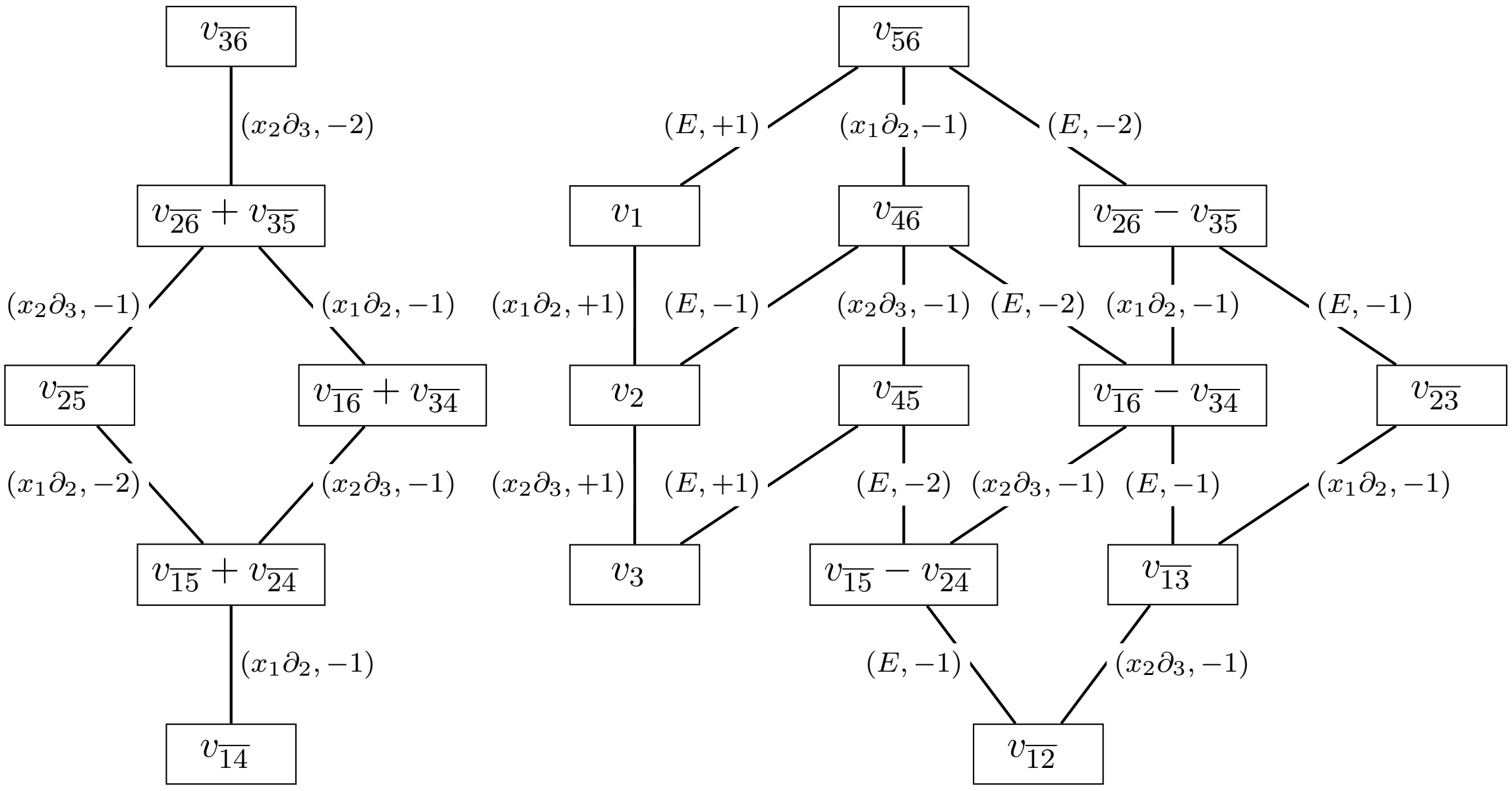

As in the degree-one case, we visualize the conditions
\[
(x_1\partial_2)w_{\mathrm{sing}}(2)=0,
\qquad
(x_2\partial_3)w_{\mathrm{sing}}(2)=0,
\qquad
Ew_{\mathrm{sing}}(2)=0
\]
in Figure~\ref{FigureDegree2}.

%\end{landscape}

\medskip
The weights of the vectors $v_i$ and $v_{\overline{ij}}$,
together with the corresponding steps if they are highest-weight vectors, are listed in
Table~\ref{weightsDegree2}.

\bigskip

\begin{table}
\begin{tabular}[]{|l|l|l|}
	\hline
	Elements& Weight& Step\\
	\hline
	\hline 
	 $v_{\overline{14}}$&	$\omega_1 =({a}-2,{b},{c}, t+2)$&$(-2,0,0,2)$ \\
	\hline
	$v_{\overline{12}}$ & $\omega_2 =({a},{b}-1,{c}-2, t+2)$& $(0,-1,-2,2)$\\
	\hline
	$v_3,v_{\overline{15}},v_{\overline{24}}$ & $\omega_3 =({a},{b}-1,{c}, t+2)$&$(0,-1,0,2)$\\
	\hline
	$v_{\overline{13}}$& $\omega_4 =({a}-1,{b}+1,{c}-2,t+2)$& $(-1,1,-2,2)$\\
	\hline
	$v_{\overline{25}}$& $\omega_5 =({a}+2,{b}-2,{c}, t+2)$& $(2,-2,0,2)$\\
	\hline 
	$v_2,v_{\overline{16}},v_{\overline{34}}$&$\omega_6 =({a}-1,{b}+1,{c},t+2)$& $(-1,1,0,2)$\\
	\hline
	$v_{\overline{45}}$&$\omega_7 =({a},{b}-1,{c}+2, t+2)$& $(0,-1,2,2)$\\
	\hline
	$ v_{\overline{23}}$&$\omega_8 =({a}+1,{b},{c}-2, t+2)$& $(1,0,-2,2)$\\
	\hline
	$ v_1,v_{\overline{26}},v_{\overline{35}}$&$\omega_9 =({a}+1,{b},{c},t+2)$& $(1,0,0,2)$\\
	\hline 
	$v_{\overline{46}}$&$\omega_{10} =({a}-1,{b}+1,{c}+2, t+2)$& $(-1,1,2,2)$\\
	\hline 
	 $v_{\overline{36}}$&$\omega_{11} =({a},{b}+2,{c},t+2)$& $(0,2,0,2)$\\
	\hline 
	$v_{\overline{56}}$&$\omega_ {12} =({a}+1,{b},{c}+2, t+2)$& $(1,0,2,2)$\\
	\hline 

\end{tabular}
\bigskip
\caption{Weights in degree two}	
\label{weightsDegree2}
\end{table}

\bigskip
As in the degree-one case, we can consider the analogous
table for the dual singular vector $w^\vee\in M(V')$.
We observe that, if the highest weight of $V$ is one of
$\omega_1,\omega_2,\omega_3,\omega_4,\omega_5,\omega_6$
or $\omega_7$, then the highest weight of $V'$ is,
respectively, $\omega_{11},\omega_{12},\omega_9,
\omega_{10},\omega_5,\omega_6$ or $\omega_8$,
where $(a,b,c,t)$ is replaced by the weight of $w^\vee$.

Therefore, up to duality, we may assume that the weight
spaces corresponding to $\omega_{11},\omega_{12},
\omega_9,\omega_{10}$ and $\omega_8$ vanish. In
particular,
\[
v_{\overline{56}}=
v_{\overline{36}}=
v_1=
v_{\overline{26}}=
v_{\overline{35}}=
v_{\overline{46}}=
v_{\overline{23}}=0.
\]
\begin{proposition}\label{sing2simple}

Let $w_{sing}(2)$ be a singular vector in a finite Verma module $M(V)$. Assume that the weight spaces corresponding to
$\omega_{11},\omega_{12},\omega_9,\omega_{10}$ and
$\omega_8$ vanish. Then, up to a scalar factor, $w_{sing}(2)$ is one of the following:
\[w[2_{AB}]
=d_1^+d_1^-\otimes x_1^{a-2}
\in M_B(a-2,0),
\]
for $a\geq 2$, or
\[
w[2_{AD}]
=d_1^+\Big(
\sum_{i=2}^3 d_i^+\otimes
\de_i\de_{z_1}\de_{z_2}
+
\sum_{j=1}^3 d_j^-\otimes
\de_j\de_{z_2}^2
\Big)
\in M_D(1,2).
\]
\end{proposition}

\begin{proof}
Applying  $x_1x_3\otimes z_+$, $d_{12}\otimes z_-$, $d_{23}\otimes z_+$,
$x_i^2\otimes z_+$ for $i=2,3$ to the vector $w_{sing}(2)$,
%$x_ix_j\otimes e_k$ and $d_{ij}\otimes e_k$ for all possible values of $i,j,k$ one obtains, among others the following system of equations
we obtain the following equations:
	\begin{align}
	\label{eq3.2} &(a+b+1)v_{\overline{25}}=0;\\
	\label{eq3.4} &(a+b)v_{\overline{45}}=0;\\
	\label{eq12.4} &(2a+4b+3c-t+12)v_{\overline{34}}=0;\\
	\label{eq15.3} & (4a+2b+3c+t+6)v_{\overline{34}}=0;\\
	\label{eq5.2} & v_2+v_{\overline{34}}=0;\\
	\label{eq9.3}& v_3+x_3 \de_2 v_{\overline{34}}=0;\\
	\label{eq4.4}& -v_3+({a}+{b}+1) v_{\overline{24}}+x_3\de_2 v_{\overline{34}}=0.
	\end{align}
We first observe that \eqref{eq3.2} implies
$v_{\overline{25}}=0$. Moreover, \eqref{eq3.4}
implies $v_{\overline{45}}=0$: indeed, if
$v_{\overline{45}}\neq 0$, then it is a highest
weight vector of weight $(a,b-1,c+2,t+2)$, so
necessarily $b\geq 1$.

Equations \eqref{eq12.4} and \eqref{eq15.3} then
imply $v_{\overline{34}}=0$. It follows from
\eqref{eq5.2}, \eqref{eq9.3} and \eqref{eq4.4},
respectively, that
\[
v_2=v_3=v_{\overline{24}}=0.
\]
We are therefore left with the five vectors
$v_{\overline{1j}}$, for $j=2,\ldots,6$.

Applying
$x_1x_2\otimes z_+$, $x_2x_3\otimes z_-$,
$d_{12}\otimes z_+$, $d_{23}\otimes z_-$,
$d_{12}\otimes z_-$, $x_3^2\otimes z_-$ and
$d_{23}\otimes z_+$ to $w_{\mathrm{sing}}(2)$,
we obtain the following equations:

	\begin{align}
	\label{eq1.1}& (a-1)v_{\overline{16}}=0\\
	%& v_2+(a-1)v_{\overline{16}}+v_{\overline{34}}=0\\
	\label{eq8.6} &  {b} v_{\overline{16}}=0\\
	% v_2 +{b} v_{\overline{16}} +v_{\overline{34}}=0\\
	\label{eq11.1}& (-2{a}-4{b}+3{c}+t-6)v_{\overline{16}}=0\\
	%-2{a}-4{b}+3{c}+t-6)v_{\overline{16}}-6v_{\overline{34}}=0\\
	\label{eq16.6} & (4{a}+2{b}-3{c}+t)v_{\overline{16}}=0\\
	\label{eq16.3} &  (4{a}+2{b}-3{c}+t)v_{\overline{13}}+6Fv_{\overline{34}}=0\\
	%& (4{a}+2{b}-3{c}+t)v_{\overline{16}}+6v_{\overline{34}}=0\\
	\label{eq16.1}& x_2 \de_1 v_{\overline{12}}+x_3 \de_1 v_{\overline{13}}-Fv_{\overline{14}}=0	\\
	\label{eq12.1} & (2{a}+4{b}+3{c}-t+6)v_{\overline{13}}+6Fv_{\overline{16}}=0\\
\label{eq10.6} &v_{\overline{15}}+x_3 \de_2 v_{\overline{16}}=0\\
		\label{eq10.3} &v_{\overline{12}} +x_3 \de_2 v_{\overline{13}}=0;\\
	\label{eq15.1} &(4{a}+2{b}+3{c}+t)v_{\overline{14}}+6x_2 \de_1 v_{\overline{15}}+6x_3 \de_1 v_{\overline{16}}=0.
	\end{align}
	
	If $v_{\overline{16}}\neq 0$, then $v_{\overline{16}}$ is a highest weight vector. We note that \eqref{eq1.1} gives ${a}=1$, \eqref{eq8.6} gives $b=0$ while \eqref{eq11.1} and \eqref{eq16.6} give ${c}=2$ and $t=2$, hence the weight 
of $v_{\overline{16}}$ is $(0,1,2,4)$. Using equation \eqref{eq10.6}, \eqref{eq15.1},
\eqref{eq12.1} and \eqref{eq10.3} we obtain the singular vector $w[2_{AD}]$.

    Now assume that $v_{\overline{16}}=0$.
Then \eqref{eq10.6} gives $v_{\overline{15}}=0$.

Moreover, applying $x_1x_2\otimes z_-$ and
$x_2x_3\otimes z_+$ to $w_{\mathrm{sing}}(2)$,
we obtain

	\begin{align}
	\label{eq2.1} &({a}-1)v_{\overline{13}}=0\\	
	\label{eq7.1} &{b}v_{\overline{14}}=0.
	\end{align}
If $v_{\overline{13}}\neq 0$, then
$a=1$ by \eqref{eq2.1}. This contradicts
\eqref{eq12.1} and \eqref{eq16.3}, since
$b\geq 0$. Hence $v_{\overline{13}}=0$,
and \eqref{eq10.3} then gives
$v_{\overline{12}}=0$.

Thus, the only possibly nonzero term is
$v_{\overline{14}}$. Equation \eqref{eq7.1}
gives $b=0$. Since $v_{\overline{14}}$ is a highest weight vector
and its $\mathfrak{sl}_2$-weight is $c$, equation \eqref{eq16.1} gives
$c=0$. Finally, \eqref{eq15.1} yields
$t=-4a$, and we obtain the singular vector
$w[2_{AB}]$.
\end{proof}

\begin{theorem}\label{degree2} The following is a complete list of singular vectors
of degree $2$ in finite Verma modules, up to scalar factors, together with
the corresponding morphisms:
\begin{itemize}\item $w[2_{AB}]=d_1^+d_1^-\otimes x_1^{a-2} \in M_B(a-2,0)$ giving the morphism $\varphi[2_{AB}]: M_A(a,0)\rightarrow M_B(a-2,0)$, for $a\geq 2$;
\item $w[2_{AD}]=d_1^+\Big(\sum_{i=2}^3 d_i^+\otimes \de_i \de_{z_1}\de_{z_2}+
	\sum_{j=1}^3 d_j^-\otimes \de_j\de_{z_2}^2\Big) \in M(0,1,2,4)$ giving the morphism  $\varphi[2_{AD}]=\varphi[1_D]\circ \varphi[1_A]: M_A(1,2)\rightarrow M_D(1,2)$; 
	\item $w[2_{CD}]=\sum _{i,j}d_i^+d_j^-\otimes \de_i \de_j \de_3^b\in M_D(b+2,0)$ giving the morphism $\varphi[2_{CD}]: M_C(b,0)\rightarrow M_D(b+2,0)$, for $b\geq 0$.
\end{itemize}
We have $w[2_{AB}]^\vee=w[2_{CD}]$, whilst $w[2_{AD}]$ is self-dual.
\end{theorem} 
\begin{proof}
The result follows from Proposition~\ref{sing2simple}
and the preceding discussion of duality. We only
remark that the morphism $\varphi[2_{AD}]$ is
self-dual and satisfies
\[
\varphi[2_{AD}]
=
\varphi[1_D]\circ\varphi[1_A].
\]
\end{proof}
\section{Degree Three}
%{\color{blue} \begin{remark}
%	If a weight $(a,b,c,t)$ is lower than $(a',b',c',t)$, then $2a+2b+c<2a'+2b'+c'$. This is often used without explicit mention.
%\end{remark}}
In this section we study singular vectors of degree $3$.
A singular vector of degree $3$ in $M(V)$ can be written
in the following form:
\[
w_{sing}(3)=\sum_{j=7}^8 n_{{j}}\otimes v_{\overline{j}}+\sum_{i=1}^3\sum_{j=1}^6 \de_i n_{{j}}\otimes v_{i\overline{j}}+\sum_{1\leq i<j<k\leq 6}n_{{i}} n_{{j}} n_{{k}} \otimes v_{\overline{ijk}},
\]
where we recall that $n_{{7}}=1\otimes z_+$, $n_{{8}}=1\otimes z_-$ and 
$n_{i}=d_i^+$ for $i=1,2,3$ and $n_{i}=d_{i-3}^-$ for $i=4,5,6$,
and we denote by $(a,b,c,t)$ the weight of $w_{sing}(3)$.
The $40$ vectors $v_{\overline j}$,
$v_{i\overline j}$ and $v_{\overline{ijk}}$,
for the appropriate values of the indices, lie in
$16$ possible weight spaces, as shown in
Table~\ref{weightDegree3}. The last column records
the step of $w_{\mathrm{sing}}(3)$ corresponding to
each possible highest weight.

\begin{table}
\begin{tabular}[]{|l|l|l|}
\hline
Element & Weight & Step\\
\hline
\hline
$v_{1\overline {6}}, v_{\overline{356}}$ & $\omega_1=({a}+1,{b}+1,{c}+1, t+3)$ & $(1,1,1,3)$ \\
\hline 
 $v_{1\overline{3}}, v_{\overline{236}}$ & $\omega_2=({a}+1,{b}+1,{c}-1, t+3)$ & $(1,1,-1,3)$\\
\hline 
 $v_{1\overline{5}}, v_{\overline{256}}$ & $\omega_3=({a}+2,{b}-1,{c}+1,t+3)$ & $(2,-1,1,3)$\\
\hline 
 $v_{1\overline{2}}, v_{\overline{235}}$ & $\omega_4=({a}+2,{b}-1,{c}-1,t+3)$ & $(2,-1,-1,3)$\\
\hline 
$v_{2\overline{6}}, v_{\overline{346}}$ & $\omega_5=({a}-1,{b}+2,{c}+1,t+3)$ & $(-1,2,1,3)$\\
\hline 
$v_{2\overline{3}}, v_{\overline{136}}$ & $\omega_6=({a}-1,{b}+2,{c}-1,t+3)$ & $(-1,2,-1,3)$\\
\hline 
$v_{\overline{456}}$ & $\omega_7=({a},{b},{c}+3,t+3)$ & $(0,0,3,3)$\\
\hline 
$v_{\overline{8}},v_{1\overline{4}},v_{2\overline{5}},v_{3\overline{6}}, v_{\overline{156}}, v_{\overline{246}},v_{\overline{345}}$ & $\omega_8=({a},{b},{c}+1,t+3)$ & $(0,0,1,3)$\\
\hline 
$v_{\overline{7}},v_{1\overline{1}},v_{2\overline{2}},v_{3\overline{3}}, v_{\overline{135}}, v_{\overline{126}},v_{\overline{234}}$ & $\omega_9=({a},{b},{c}-1,t+3)$ & $(0,0,-1,3)$\\
\hline 
$v_{\overline{123}}$ & $\omega_{10}=({a},{b},{c}-3,t+3)$ & $(0,0,-3,3)$\\
\hline 
 $v_{3\overline{5}}, v_{\overline{245}}$ & $\omega_{11}=({a}+1,{b}-2,{c}+1,t+3)$ & $(1,-2,1,3)$\\
\hline 
 $v_{3\overline{2}}, v_{\overline{125}}$ & $\omega_{12}=({a}+1,{b}-2,{c}-1,t+3)$ & $(1,-2,-1,3)$\\
\hline 
$v_{2\overline{4}}, v_{\overline{146}}$ & $\omega_{13}=({a}-2,{b}+1,{c}+1,t+3)$ & $(-2,1,1,3)$\\
\hline 
 $v_{2\overline{1}}, v_{\overline{134}}$ & $\omega_{14}=({a}-2,{b}+1,{c}-1,t+3)$ & $(-2,1,-1,3)$\\
\hline 
 $v_{3\overline{4}}, v_{\overline{145}}$ & $\omega_{15}=({a}-1,{b}-1,{c}+1,t+3)$ & $(-1,-1,1,3)$\\
\hline 
$v_{3\overline{1}}, v_{\overline{124}}$ & $\omega_{16}=({a}-1,{b}-1,{c}-1,t+3)$ & $(-1,-1,-1,3)$\\
\hline 
\end{tabular}
\bigskip
\caption{Weights in degree three.}
\label{weightDegree3}
\end{table}
The singularity conditions obtained by applying
$x_1^2\otimes z_+$ to $w_{\mathrm{sing}}(3)$ imply
 $v_{1\overline 6}=v_{1\overline 5}=v_{\overline{256}}=v_{\overline{356}}=0$. Applying $x_1^2\otimes z_-$ one gets $v_{1\overline 3}=v_{1\overline 2}=v_{\overline{236}}=v_{\overline{235}}=0$, applying $x_1x_2\otimes z_-$ we get $v_{2\overline 3}=0$, applying $d_{12}\otimes z_-$ we get $v_{\overline {346}}=0$ and finally, applying $E$ we also get $v_{2\overline 6}=0$.
Therefore we are left with 29 possible nonzero summands, listed in Table
\ref{grado3reduced}.

\begin{table}
\begin{tabular}[]{|c|c|c|}
	\hline
Element	& Weight & Step\\
	\hline
	\hline 
	 $ v_{\overline{136}}$ & $\omega_6=({a}-1,{b}+2,{c}-1,t+3)$ & $(-1,2,-1,3)$\\
	\hline 
 $v_{\overline{456}}$ & 	$\omega_7=({a},{b},{c}+3,t+3)$ & $(0,0,3,3)$\\
	\hline 
 $v_{\overline{8}},v_{1\overline{4}},v_{2\overline{5}},v_{3\overline{6}}, v_{\overline{156}}, v_{\overline{246}},v_{\overline{345}}$ & 	$\omega_8=({a},{b},{c}+1,t+3)$ & $(0,0,1,3)$\\
	\hline 
$v_{\overline{7}},v_{1\overline{1}},v_{2\overline{2}},v_{3\overline{3}}, v_{\overline{135}}, v_{\overline{126}},v_{\overline{234}}$ & 	$\omega_9=({a},{b},{c}-1,t+3)$ & $(0,0,-1,3)$\\
	\hline 
$v_{\overline{123}}$ & 	$\omega_{10}=({a},{b},{c}-3,t+3)$ & $(0,0,-3,3)$\\
	\hline 
 $v_{3\overline{5}}, v_{\overline{245}}$ & 	$\omega_{11}=({a}+1,{b}-2,{c}+1,t+3)$ & $(1,-2,1,3)$\\
	\hline 
 $v_{3\overline{2}}, v_{\overline{125}}$ & 	$\omega_{12}=({a}+1,{b}-2,{c}-1,t+3)$ & $(1,-2,-1,3)$\\
	\hline 
 $v_{2\overline{4}}, v_{\overline{146}}$ & 	$\omega_{13}=({a}-2,{b}+1,{c}+1,t+3)$ & $(-2,1,1,3)$\\
	\hline 
$v_{2\overline{1}}, v_{\overline{134}}$ & 	$\omega_{14}=({a}-2,{b}+1,{c}-1,t+3)$ & $(-2,1,-1,3)$\\
	\hline 
	 $v_{3\overline{4}}, v_{\overline{145}}$ & $\omega_{15}=({a}-1,{b}-1,{c}+1,t+3)$ &
	 $(-1,-1,1,3)$\\
	\hline 
 $v_{3\overline{1}}, v_{\overline{124}}$& 	$\omega_{16}=({a}-1,{b}-1,{c}-1,t+3)$ & 
  $(-1,-1,-1,3)$\\
	\hline 
\end{tabular}
\bigskip
\caption{Possible nonzero components of $w_{\mathrm{sing}}(3)$.}
\label{grado3reduced}
\end{table}

\bigskip

The existence of duality between singular vectors allows us to make the following crucial observation.
\begin{lemma}
	No weight among $\omega_{11}$, $\omega_{12}$,  $\omega_{14}$, $\omega_{15}$ and $\omega_{16}$ can be the highest weight in $V$.
\end{lemma}
	\begin{proof}
	The result follows, by duality, from Remark \ref{rem_dual} and Table \ref{grado3reduced}.
%		Recall the duality theorem: if there exists a morphism
%		\[
%		M(a,b,c,t)\rightarrow M(a',b',c',t+3)
%		\]
%		then there exists the dual morphism,
%		\[
%		M(b',a',c', -t+3)\rightarrow M(b,a,c,-t+6)
%		\]
%	Assume that $(a-1,b-1,c-1,t+3)$ is the highest weight so that we have a morphism
%	\[
%	M(a,b,c,t)\rightarrow M(a-1,b-1,c-1,t+3)
%	\]
%	then the dual morphism is
%	\[
%	M(b-1,a-1,c-1,-t+3)\rightarrow M(b,a,c,-t+6)
%	\]
%	If we denote by $(a^*,b^*,c^*,t^*)=(b-1,a-1,c-1,-t+3)$ the weight of the dual singular vector we may rewrite the dual morphism as
%	\[
%	M(a^*,b^*,c^*,t^*)\rightarrow (a^*+1,b^*+1,c^*+1,t^*+3)
%	\]
%	and this is a contradiction since we do not have terms of weight $(a^*+1,b^*+1,c^*+1,t^*+3)$ in the table 2 (the only possible such terms appearing in table 1 vanish). The same argument applies to the other weights. 
	\end{proof}

\begin{lemma}
	We have $v_{\overline{136}}=0$. In particular $\omega_6$ is not the highest weight in $V$ and, by duality, $\omega_{13}$ is not the highest weight either.
\end{lemma}
\begin{proof}
Applying $x_1x_2\otimes z_+$, $x_2x_3\otimes z_-$, $d_{12}\otimes z_+$,  $d_{12}\otimes z_-$ and $d_{23}\otimes z_-$ we obtain the following equations
\[
\begin{cases}
	(a-1)v_{\overline{136}}=0\\
	bv_{\overline{136}}=0\\
	(2a + 4b - 3c - t + 12) v_{\overline{136}}=0\\
	Fv_{\overline{136}}=0\\
	(4a + 2b - 3c + t) v_{\overline{136}}=0
\end{cases}
\]
which are compatible if and only if $v_{\overline{136}}=0$, since the condition $Fv_{\overline{136}}=0$ implies $c=1$ if $v_{\overline{136}}\neq 0$.
\end{proof}

\begin{lemma}
The weight	$\omega_9$ is not the highest weight and by duality $\omega_8$ is not the highest weight either.
\end{lemma}
\begin{proof}
	Applying $x_1^2\otimes z_-$, $x_1x_2\otimes z_-$, $x_1x_3\otimes z_-$,  $x_2x_3\otimes z_+$, $x_3^2\otimes z_-$, $d_{12}\otimes z_+$, $d_{12}\otimes z_-$ 
	%(two equations) 
	and $x_1x_3\otimes z_+$  we obtain the following equations:
	\begin{align}
		&\label{eqdeg31} v_{1\overline 1}=-v_{\overline{234}};\\ 
		&\label{eqdeg32}v_{2\overline 2}=-({a}+1)v_{\overline{234}};\\
		&\label{eqdeg33}v_{3\overline 3}=-({a}+{b}+1)v_{\overline{234}}\\
		&\label{eqdeg34}{a}v_{3\overline 3}=0\\
		&\label{eqdeg35}v_{2\overline 2}-	v_{3\overline 3}+{b}v_{\overline{234}}=0\\
		&\label{eqdeg36}v_{1\overline 1}+(a-1)v_{\overline{135}}=0\\
		&\label{eqdeg37}v_{3\overline 3}-v_{\overline{126}}-v_{\overline{135}}+v_{\overline{234}}=0\\
		&\label{eqdeg38}(2{a} + 4{b} - 3c - t + 6)v_{\overline{126}}+6v_{\overline{135}}-6v_{\overline{234}}=0\\
		&\label{eqdeg39}6v_{\overline 7}+6v_{1\overline 1}+( -2{a}-4{b}-3c+t-6)v_{3\overline{3}}+6v_{2\overline{2}}=0\\
		&\label{eqdeg310}(2{a}+4{b}+3{c}-t+6)v_{\overline{234}}=0\\
		&\label{eqdeg311} v_{1\overline 1}-v_{3\overline 3}-({a}+{b}+1)v_{\overline{135}}-v_{\overline{234}}=0		
	\end{align}

Assume that $v_{\overline{234}}\neq 0$.
By \eqref{eqdeg31}--\eqref{eqdeg33}, the vectors
$v_{1\overline 1}$, $v_{2\overline 2}$ and
$v_{3\overline 3}$ are also nonzero. Equation
\eqref{eqdeg34} gives $a=0$, while
\eqref{eqdeg35} then yields $b=0$.
Equations \eqref{eqdeg36} and \eqref{eqdeg37}
give
\[
v_{\overline{135}}=-v_{\overline{234}},
\qquad
v_{\overline{126}}=v_{\overline{234}}.
\]
Now \eqref{eqdeg38} and \eqref{eqdeg310}
respectively imply
\[
t=-3c-6,\qquad t=3c+6,
\]
which is impossible since $c\geq 0$.
Therefore $v_{\overline{234}}=0$. It follows from \eqref{eqdeg31}--\eqref{eqdeg33}
that
\[
v_{1\overline 1}=v_{2\overline 2}
=v_{3\overline 3}=0.
\]
Equation \eqref{eqdeg311} then gives
$v_{\overline{135}}=0$, while
\eqref{eqdeg37} yields $v_{\overline{126}}=0$.
Finally, \eqref{eqdeg39} implies
$v_{\overline 7}=0$.

\end{proof}
Recall that we let $D_3=d_1^+d_2^+d_3^+$, $D_2=d_1^-d_2^+d_3^++d_1^+d_2^-d_3^++d_1^+d_2^+d_3^-$, $D_1=d_1^+d_2^-d_3^-+d_1^-d_2^+d_3^-+d_1^-d_2^-d_3^+$ and $D_0=d_1^-d_2^-d_3^-$.
\begin{theorem}\label{ClassDeg3}The following is a complete list, up to scalar factors,  of singular vectors of degree 3 and corresponding morphisms of finite Verma modules:
	\begin{itemize}
		\item $w[3_{AC}]=D_3\otimes z_1^{c-3}\in M_A(0,c-3)$ for all $c\geq 3$, inducing the morphism $\varphi[3_{AC}]: M_A(0,c)\rightarrow M_C(0,c-3)$;\vspace{2mm}
		\item $w[3_{BD}]=\sum_{k=0}^3 D_k\otimes \de_+^k \de_-^{c+3-k}\in M_D(0,c+3)$ 
        %\begin{align*}w[3_{BD}]=& d_1^-d_2^-d_3^-\otimes z_1^c-(d_1^+d_2^-d_3^-+d_1^-d_2^+d_3^-+d_1^-d_2^-d_3^+)\otimes z_1^{c-1}z_2\\&+(d_1^-d_2^+d_3^++d_1^+d_2^-d_3^++d_1^+d_2^+d_3^-)\otimes z_1^{c-2}z_2^2-d_1^+d_2^+d_3^+\otimes z_1^{c-3}z_2^3  \in M(0,0,c,-3c+6) \end{align*} 
        for all $c\geq 0$ inducing the morphism
		$\varphi[3_{BD}]: M_B(0,c)\rightarrow M_D(0,c+3)$.
	\end{itemize} 

\end{theorem}
\begin{proof}
Denote, as above, the weight of the singular vector by $(a,b,c,t)$.  By the previous results the only two possible highest weights in $V$ are $({a},{b},{c}-3,{t+3})$ and $({a},{b},{c}+3,{t+3})$ and, by duality, there is a correspondence between the two cases.
If the highest weight is $({a},{b},{c}-3,{t+3})$ then necessarily the singular vector is of the form $n_{1}n_{2}n_{3}\otimes v_{\overline{123}}=d_1^+d_2^+d_3^+\otimes v_{\overline{123}}$. Applying $x_1x_3\otimes z_-$ we get the equations $(a+b)v_{\overline{123}}=0$ forcing $a=b=0$ and applying $d_{12}\otimes z_+$ we also deduce $t=3c$ and we get the vector $w[3_{AC}]$, which is indeed singular. 

The result follows by duality, since a direct computation shows that $w[3_{BD}]$ is also singular.
\end{proof}

\section{Degree 4}\label{S7}
In this section we give a sketch of the proof of the
classification of singular vectors of degree $4$.
A singular vector of degree $4$ and weight $(a,b,c,t)$
in $M(V)$ can be written in the following form:

\begin{align*}
w_{sing}(4)&=\sum_{1\leq i<j<h<k\leq 6} n_{i}n_{j}n_{h}n_{k}\otimes v_{\overline{ijhk}}+\sum_{i=1}^3\sum_{1\leq h<k\leq 6} \de_i n_{h}n_{k}\otimes v_{i\overline{hk}}+\sum _{1\leq i\leq j\leq 3}\de_i\de_j\otimes v_{ij}\\
&+\sum_{i=7}^8\sum_{j=1}^6 n_{i}n_{j}\otimes v_{\overline{ij}}\in M(V).
\end{align*}
Here, as in degrees $2$ and $3$, we set
$n_7=d^+$, $n_8=d^-$, and
$n_i=d_i^+$ for $i=1,2,3$, while
$n_i=d_{i-3}^-$ for $i=4,5,6$.
Thus, $w_{\mathrm{sing}}(4)$ consists of $78$ terms.

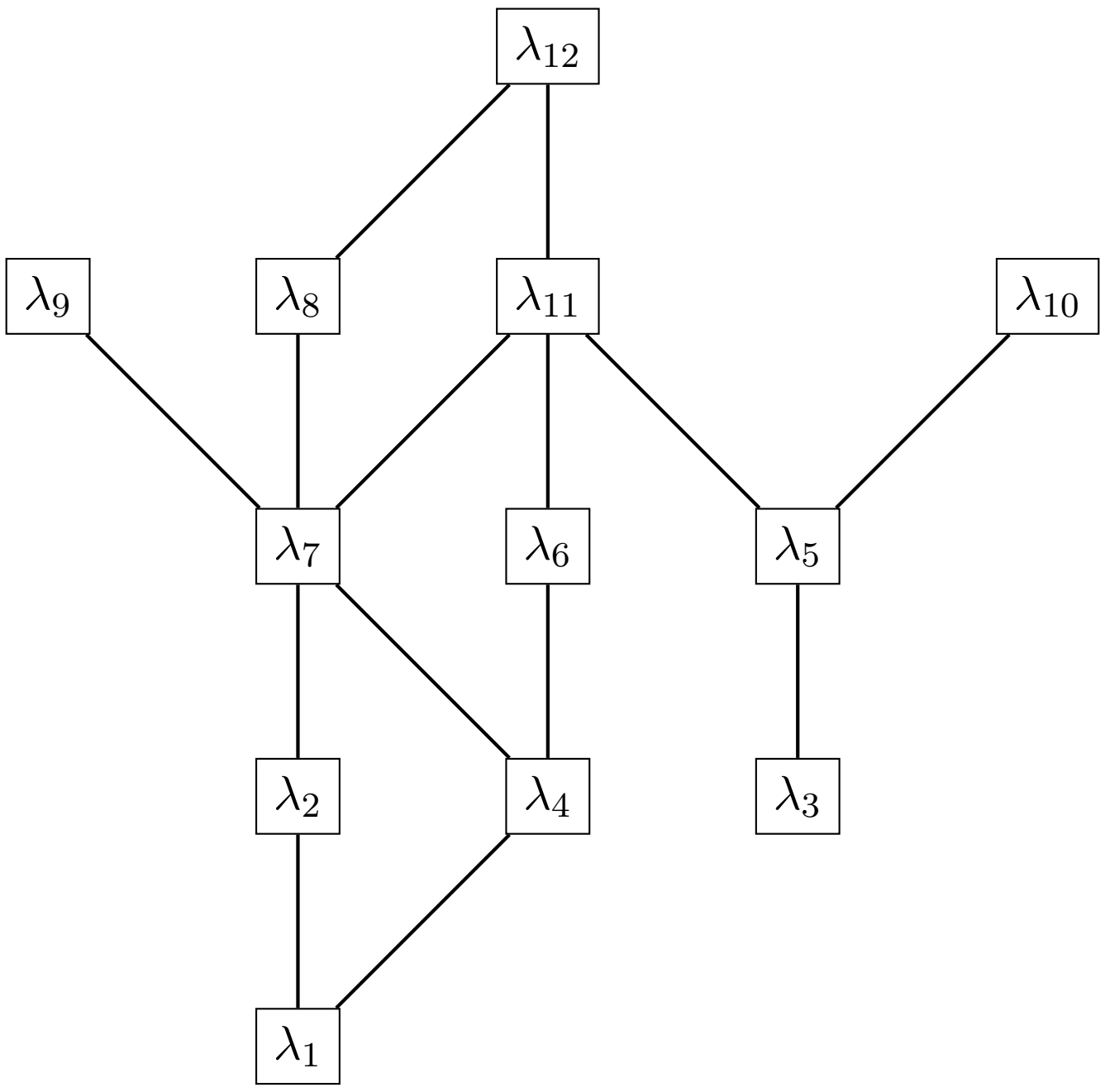
\begin{figure}[h]
\begin{tikzpicture}[
	box/.style={draw,align=center,font=\small},
	edge/.style={draw,thick},
	]
	
	\node[box] (A1) at (0,0) {$\lambda_1$};
	\node[box] (A2) at (0,2){$\lambda_2$};
	\node[box] (A3) at (4,2){$\lambda_3$};
	\node[box] (A4) at (2,2){$\lambda_4$};
	\node[box] (A5) at (4,4){$\lambda_5$};
	\node[box] (A6) at (2,4){$\lambda_6$};
	\node[box] (A7) at (0,4){$\lambda_7$};
	\node[box] (A8) at (0,6){$\lambda_8$};
	\node[box] (A9) at (-2,6){$\lambda_9$};
	\node[box] (A10) at (6,6){$\lambda_{10}$};
	\node[box] (A11) at (2,6){$\lambda_{11}$};
	\node[box] (A12) at (2,8){$\lambda_{12}$};
	
	\draw[edge] (A1) --  (A2);
	\draw[edge] (A1) --  (A4);
	\draw[edge] (A2) --  (A7);
	\draw[edge] (A4) --  (A7);
	\draw[edge] (A4) --  (A6);
	\draw[edge] (A3) --  (A5);
	\draw[edge] (A7) --  (A9);
	\draw[edge] (A7) --  (A8);
	\draw[edge] (A7) --  (A11);
	\draw[edge] (A6) --  (A11);
	\draw[edge] (A5) --  (A11);
	\draw[edge] (A5) --  (A10);
	\draw[edge] (A8) --  (A12);
	\draw[edge] (A11) --  (A12);
\end{tikzpicture}
\caption{The dominance order on weights $\lambda_1$--$\lambda_{12}$ in degree 4.}
\label{dominanceDegree4}
\end{figure}

The basic linear system in degree $4$ immediately
implies that $36$ of these terms vanish. The weights
of the remaining terms are listed in
Table~\ref{weightsDegree4}, together with the
corresponding steps of $w_{\mathrm{sing}}(4)$. Figure~\ref{dominanceDegree4} illustrates the
dominance order among these weights. 

\bigskip

\begin{table}
\begin{tabular}[]{|m{5cm}|m{5cm}|l|}
	\hline
 Element & 	Weight & Step\\
	\hline
	\hline
 $v_{3\overline {14}}$ & 	$\lambda_1=({a}-2,{b}-1,{c},t+4)$ & $(-2,-1,0,4)$ \\
	\hline 
 $v_{2\overline{14}}$ & 	$\lambda_2=({a}-3,{b}+1,{c},t+4)$& $(-3,1,0,4)$ \\
	\hline 
	 $v_{\overline{17}}, v_{2\overline{12}},v_{3\overline{13}},v_{\overline{1234}}$ &
	 $\lambda_3=({a}-1,{b},{c}-2,t+4)$& $(-1,0,-2,4)$ \\
	\hline 
$v_{33}$, $v_{3\overline{15}}$, $v_{3\overline{24}}$, $v_{\overline{1245}}$
& 	$\lambda_4=({a},{b}-2,{c},t+4)$& $(0,-2,0,4)$ \\
	\hline 
 $v_{\overline{27}}, v_{\overline{1235}}$ & 	$\lambda_5=({a}+1,{b}-1,{c}-2,t+4)$&$(1,-1,-2,4)$\\
	\hline 
 $v_{3\overline{25}}$ & 	$\lambda_6=({a}+2,{b}-3,{c},t+4)$& $(2,-3,0,4)$ \\
	\hline 
$v_{\overline{18}}$, $v_{\overline{47}}$, $v_{1\overline{14}}$, $v_{2\overline{15}}$, $v_{2\overline{24}}$, 
$v_{3\overline{16}}$, $v_{3\overline{34}}$, $v_{\overline{1246}}$, $v_{\overline{1345}}$
 & 	$\lambda_7=({a}-1,{b},{c},t+4)$& $(-1,0,0,4)$ \\
	\hline 
 $v_{2\overline{16}}$, $v_{\overline{1346}}$ & 	$\lambda_8=({a}-2,{b}+2,{c},t+4)$ &
 $(-2,2,0,4)$\\
	\hline 
$v_{\overline{1456}}$ &  	$\lambda_9=({a}-1,{b},{c}+2,t+4)$ & $(-1,0,2,4)$ \\
	\hline 
 $v_{\overline{1236}}$& 	$\lambda_{10}=({a},{b}+1,{c}-2,t+4)$ & $(0,1,-2,4)$ \\
	\hline 
 $v_{\overline{28}}$,$v_{\overline{57}}$, $v_{1\overline{15}}$, $v_{1\overline{24}}$, $v_{2\overline{25}}$,
  $v_{3\overline{26}}$, $v_{3\overline{35}}$, $v_{\overline{1256}}$, $v_{\overline{2345}}$ & 	$\lambda_{11}=({a}+1,{b}-1,{c},t+4)$& $(1,-1,0,4)$\\
	\hline 
 $v_{\overline{38}}$, $v_{\overline{67}}$, $v_{1\overline{16}}$, $v_{2\overline{26}}$, $v_{3\overline{36}}$, $v_{\overline{1356}}$, $v_{\overline{2346}}$ &
 	$\lambda_{12}=({a},{b}+1,{c},t+4)$& $(0,1,0,4)$ \\
	
	\hline 
\end{tabular}
\bigskip
\caption{\label{weightsDegree4} Weights in degree four.}

\end{table}

\bigskip

\begin{lemma}\label{poss_weights4}
	The highest weight of $V$ is one of the following: $\lambda_7,\lambda_8, \lambda_9, \lambda_{10},\lambda_{11},\lambda_{12}$.
\end{lemma}
\begin{proof}
The result follows from Remark~\ref{rem_dual} and
Table~\ref{weightsDegree4}. More precisely, duality
exchanges $\lambda_7$ with $\lambda_{12}$ and
$\lambda_9$ with $\lambda_{10}$. Moreover,
$\lambda_8$ and $\lambda_{11}$ are self-dual, whilst the duals of the other weights do not appear in the list.
%	as in the case of degree 3 by duality reasons. For example, if $\lambda_1=(a-2,b-1,c)$ is the highest weight, then considering the dual singular vector we should have that the highest weight is $(a+1,b+2,c)$ which does not appear in the table. The only possibilities are indeed $\lambda_8$ and $\lambda_{11}$, which are self-dual, and the two pairs of dual weights $\lambda_9$-$\lambda_{10}$ and $\lambda_7$-$\lambda_{12}$.
\end{proof}

\begin{lemma}The highest weight of $V$ is either $\lambda_9$ or $\lambda_{10}$. 
\end{lemma}
\begin{proof}
Since $\lambda_7$ and $\lambda_{12}$ are dual to
each other, by duality we may assume that
$\lambda_{12}$ is not the highest weight. Hence
all vectors of weight $\lambda_{12}$ must vanish.

If a nonzero vector of weight $\lambda_9$ or
$\lambda_{10}$ occurs, the result follows.
Otherwise, we may assume that all vectors of these
two weights vanish as well.

Applying $x_1x_3\otimes z_\epsilon$, for
$\epsilon=+,-$, $x_3x_1\otimes z_+$ and
$x_3^2\otimes z_-$ to $w_{\mathrm{sing}}(4)$,
we obtain that all vectors of weights
$\lambda_5$, $\lambda_6$ and $\lambda_{11}$ vanish.
In particular, $\lambda_{11}$ cannot be the
highest weight.

We have to show that all terms of weight $\lambda_7$ or $\lambda_8$ vanish.

Condition $(x_2^2\otimes z_-)w_{sing}(4)=0$ gives $v_{2\overline{16}}=-v_{\overline{1346}}$. If we assume, by contradiction, that $v_{2\overline{16}}=-v_{\overline{1346}}$ are nonzero, conditions
$(x_1x_2\otimes z_+)w_{sing}(4)=0$ and $(x_2x_3\otimes z_-)w_{sing}(4)=0$ provide $a=2$ and $b=0$, respectively. Then applying $d_{23}\otimes z_-$ to $w_{sing}(4)$ one also has $Fv_{2\overline{16}}=0$ which forces $c=0$. Finally, applying $d_{23}\otimes z_+$ and $d_{12}\otimes z_+$ to $w_{sing}(4)$  one gets a contradiction.
%One also has $(8 + t)v_{2\overline{16}}=0$ and $(t - 16)v_{2\overline{16}} -6v_{\overline{1346}}=0$ achieving a contradiction. So we have $v_{2\overline{16}}=v_{\overline{1346}}=0$. In particular the two terms of weight $\lambda_8$ vanish and so $\lambda_8$ is not the highest weight.

Looking at the equations given by $(x_2^2\otimes z_-)w_{sing}(4)=0$ and 
$(x_1x_3\otimes z_-)w_{sing}(4)=0$, $(x_2x_3\otimes z_+)w_{sing}(4)=0$, we have $v_{2\overline{24}}=v_{3\overline{34}}=0$ and $v_{33}=v_{3\overline{24}}=0$; moreover one can notice, using conditions $(x_2x_3\otimes z_-)w_{sing}(4)=0$ and $(x_2x_3\otimes z_+)w_{sing}(4)=0$,  that if $b\neq 0$ then all terms vanish, so we can assume $b=0$. The same equations imply $v_{2\overline{15}}=v_{\overline{1246}}=v_{1\overline{14}}=v_{3\overline{16}}=0$. At this stage only two  possible nonzero terms of weight $\lambda_7$ (with $b=0$) are left: $v_{\overline {18}}$ and $v_{\overline{47}}$. But one can easily verify that they also vanish by application of $d_{12}\otimes z_+$ and $d_{12}\otimes z_-$ respectively. In particular $\lambda_7$ is not the highest weight.
\end{proof}
\begin{theorem}\label{ClassDeg4}The following is a complete list, up to a scalar factor, of singular vectors of degree 4 and corresponding morphisms of Verma modules.
\begin{itemize}
    \item $w[4_{AD}]=D_3 \sum_{i=1}^3 d_i^-\otimes \de_i \in M_D(1,0)$  inducing the morphism \\ $\varphi[4_{AD}]:M_A(0,2)\rightarrow M_D(1,0)$;\vspace{2mm}
    \item  $w'[4_{AD}]= d_1^- \sum_{k=0}^2 D_{k+1}\otimes \de_+^k \de_- ^{2-k} \in M_D(0,2)$ inducing the morphism \\ $\varphi'[4_{AD}]:M_A(1,0)\rightarrow M_D(0,2)$.
\end{itemize}

\end{theorem}
\begin{proof}
Assume that the highest weight is $\lambda_{10}$, so that $v_{\overline{1236}}$ is a highest weight vector. We observe that the only weights in $\lambda_1,\ldots,\lambda_{12}$ which are strictly lower than $\lambda_{10}$ are $\lambda_3$ and $\lambda_5$ and therefore all vectors of weight $\lambda_{i}$, for all $i\neq 3,5,10$, vanish.
The equations given by $x_ix_j\otimes z_\epsilon$ and $d_{ij}\otimes z_\epsilon$, for $\epsilon=+,-$, provide $(a+b)v_{\overline{1236}}=0$ which implies $a=b=0$. We also obtain $v_{2\overline{12}}=v_{3\overline{13}}=v_{\overline{17}}=v_{\overline{27}}=0$. The two equations $(3c+t-12)v_{\overline{1236}}=0$ and $(-3c+t)v_{\overline{1236}}=0$ force $t=6$ and $c=2$. 
Finally, the relations
\[
x_3\de_2v_{\overline{1236}}
=-v_{\overline{1235}},
\qquad
x_2\de_1v_{\overline{1235}}
=-v_{\overline{1234}}
\]
show that, up to a scalar factor, the singular
vector is $w[4_{AD}]$, which is indeed singular.
The result follows by duality, since a direct
computation shows that $w'[4_{AD}]$ is also singular.
\end{proof}

\section{Degree 6}
 We write a singular vector of weight $(a,b,c,t)$ and degree 6 in $M(V)$ in the following form:
\begin{align*}
w_{sing}(6)&=\sum_{\substack{1\leq i_1<\cdots<i_k\leq 8\\deg(n_{i_1}\cdots n_{i_k})=6}} n_{i_1}\cdots n_{i_k}\otimes v_{\overline{i_1\cdots i_k}}+\sum_{i=1}^3\sum_{\substack{1\leq i_1<\cdots<i_k\leq 8\\deg(n_{i_1}\cdots n_{i_k})=4}} \de_i n_{i_1}\cdots n_{i_k}\otimes v_{i\,\overline{i_1\cdots i_k}}\\
&+ \sum_{1\leq i\leq j\leq 3}\sum_{1\leq h<k\leq 6} \de_i\de_j n_{h}n_{k}\otimes v_{ij\,\overline{hk}}+\sum _{1\leq i\leq j\leq k\leq 3}\de_i\de_j\de_k\otimes v_{ijk}\in M(V).
\end{align*}

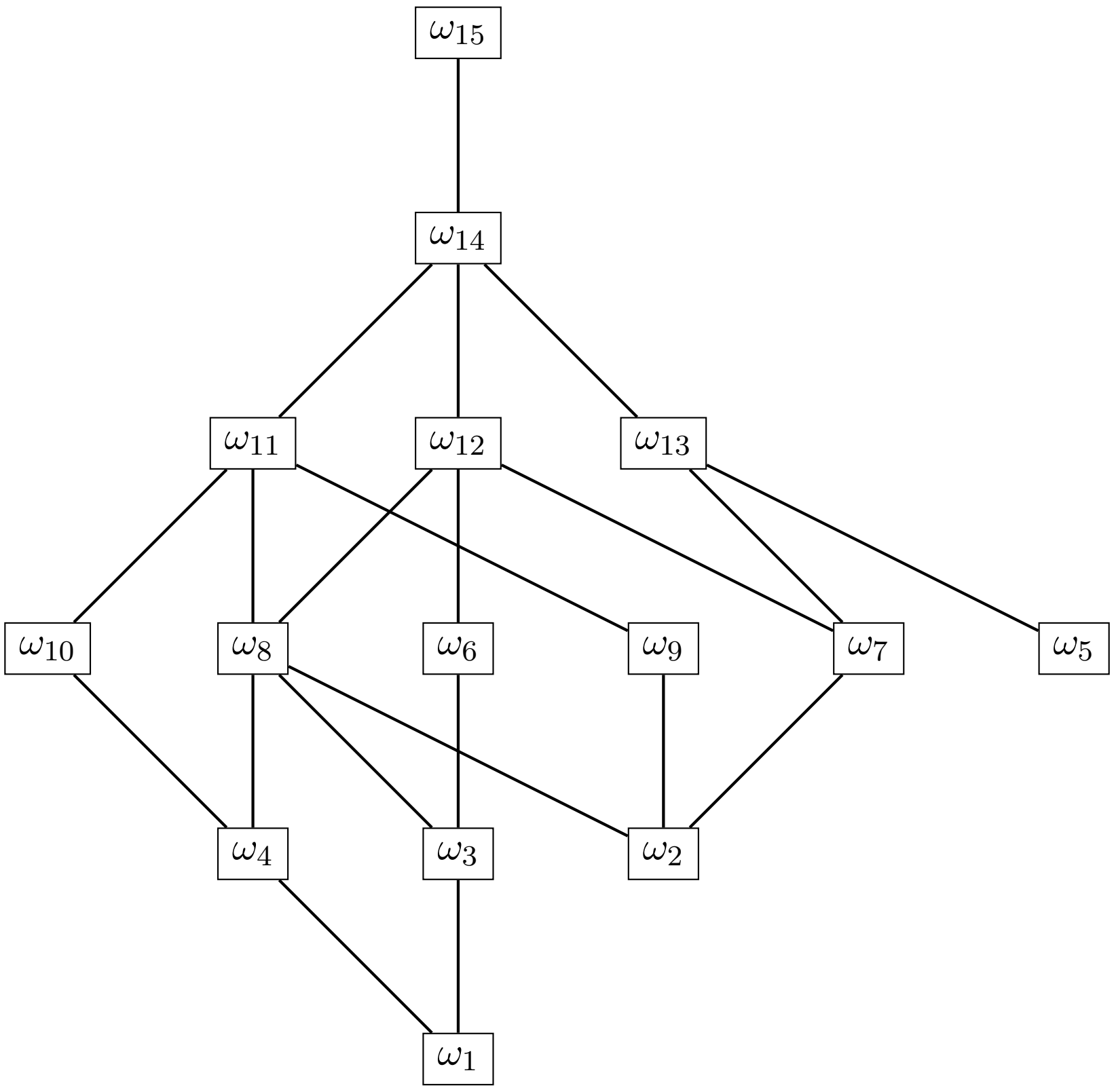
\begin{figure}[h]
	\begin{tikzpicture}[
		box/.style={draw,align=center,font=\small},
		edge/.style={draw,thick},
		]
		
		\node[box] (A1) at (0,0) {$\omega_1$};
		\node[box] (A2) at (2,2){$\omega_2$};
		\node[box] (A3) at (0,2){$\omega_3$};
		\node[box] (A4) at (-2,2){$\omega_4$};
		\node[box] (A5) at (6,4){$\omega_5$};
		\node[box] (A6) at (0,4){$\omega_6$};
		\node[box] (A7) at (4,4){$\omega_7$};
		\node[box] (A8) at (-2,4){$\omega_8$};
		\node[box] (A9) at (2,4){$\omega_9$};
		\node[box] (A10) at (-4,4){$\omega_{10}$};
		\node[box] (A11) at (-2,6){$\omega_{11}$};
		\node[box] (A12) at (0,6){$\omega_{12}$};
		\node[box] (A13) at (2,6){$\omega_{13}$};
		\node[box] (A14) at (0,8){$\omega_{14}$};
		\node[box] (A15) at (0,10){$\omega_{15}$};
		
		\draw[edge] (A1) --  (A4);
		\draw[edge] (A1) --  (A3);
		\draw[edge] (A4) --  (A10);
		\draw[edge] (A4) --  (A8);
		\draw[edge] (A3) --  (A6);
		\draw[edge] (A3) --  (A8);
		\draw[edge] (A2) --  (A8);
		\draw[edge] (A2) --  (A7);
		\draw[edge] (A2) --  (A9);
		\draw[edge] (A10) --  (A11);
		\draw[edge] (A6) --  (A12);
		\draw[edge] (A8) --  (A11);
		\draw[edge] (A8) --  (A12);
		\draw[edge] (A7) --  (A12);
		\draw[edge] (A7) --  (A13);
		\draw[edge] (A9) --  (A11);
		\draw[edge] (A5) --  (A13);
		\draw[edge] (A11) --  (A14);
		\draw[edge] (A12) --  (A14);
		\draw[edge] (A13) --  (A14);
		\draw[edge] (A14) --  (A15);

	\end{tikzpicture}
	
	\caption{\label{orderDegree6} The dominance order on the weights $\omega_1,\ldots,\omega_{15}$ in degree $6$.}
\end{figure}

In degree 6 the basic linear homogeneous system reduces the generic form of a singular vector to a sum of 89 terms distributed in 15 weight spaces as shown in Table \ref{weightsDegree6}. The dominance order of these weights is shown in Figure
	\ref{orderDegree6}.

\begin{table}
\begin{tabular}[]{|p{6cm}|m{6cm}|l|}
	\hline
 Element &	Weight\\
	\hline
	\hline	
	$v_{33\overline {14}}$ & $\omega_1=(a-2,b-2,c,t+6)$\\
	\hline 
	 $v_{3\overline{17}}$, $v_{\overline{1247}}$ & $\omega_2=(a-1,b-1,c-2,t+6)$\\
	\hline 
	 $v_{23\overline{14}}$ & $\omega_3=(a-3,b,c,t+6)$\\
	\hline 
	$v_{333}$, $v_{33\overline{15}}$, $v_{33\overline{24}}$, $v_{\overline{1245}}$ & $\omega_4=(a,b-3,c,t+6)$\\
	\hline 
	$ v_{\overline{1237}}$ & $\omega_5=(a,b,c-4,t+6)$\\
	\hline 
	$v_{22\overline{14}}$ & $\omega_6=(a-4,b+2,c,t+6)$\\
	\hline 
	$v_{2\overline{17}}$, $v_{\overline{1347}}$ & $\omega_7=(a-2,b+1,c-2,t+6)$\\
	\hline 
	$v_{233}$, $v_{3\overline{18}}$, $v_{3\overline{47}}$,$v_{13\overline{14}}$,$v_{23\overline{15}}$,$v_{23\overline{24}}$, &\\ $v_{33\overline{16}}$,$v_{33\overline{34}}$,$v_{\overline{1248}}$,$v_{\overline{1457}}$, $v_{2\overline{1245}}$,$v_{3\overline{1246}}$, $v_{3\overline{1345}}$ & $\omega_8=(a-1,b-1,c,t+6)$\\
	\hline 
	$v_{\overline{1257}}$ & $\omega_9=(a+1,b-2,c-2,t+6)$ \\
	\hline 
	$v_{33\overline{25}}$ & $\omega_{10}=(a+2,b-4,c,t+6)$\\
	\hline 
	$v_{3\overline{57}}$,$v_{13\overline{15}}$,$v_{23\overline{25}}$,$v_{33\overline{26}}$,$v_{33\overline{35}}$,$v_{\overline{1258}}$, &\\$v_{\overline{2457}}$,$v_{1\overline{1245}}$,$v_{3\overline{1256}}$,$v_{3\overline{2345}}$ & $\omega_{11}=(a+1,b-2,c,t+6)$\\
	\hline 
	$v_{2\overline{18}}$,$v_{2\overline{47}}$,$v_{12\overline{14}}$,$v_{22\overline{15}}$,$v_{23\overline{16}}$,$v_{\overline{1348}}$,$v_{\overline{1467}}$,&\\$v_{2\overline{1246}}$,$v_{2\overline{1345}}$,$v_{3\overline{1346}}$ & $\omega_{12}=(a-2,b+1,c,t+6)$\\
	\hline
	$v_{1\overline{17}}$,$v_{2\overline{27}}$,$v_{3\overline{37}}$,$v_{\overline{1238}}$,$v_{\overline{1267}}$,$v_{\overline{1357}}$,$v_{\overline{2347}}$, &\\
	$v_{1\overline{1234}}$,$v_{2\overline{1235}}$,$v_{3\overline{1236}}$ & $\omega_{13}=(a,b,c-2,t+6)$\\
	\hline
	$v_{\overline{78}}$,$v_{1\overline{18}}$,$v_{1\overline{47}}$,$v_{2\overline{28}}$,$v_{2\overline{57}}$,$v_{3\overline{38}}$,$v_{3\overline{67}}$,$v_{11\overline{14}}$,&\\$v_{12\overline{15}}$,$v_{12\overline{24}}$,$v_{13\overline{16}}$, $v_{13\overline{34}}$,$v_{22\overline{25}}$,$v_{23\overline{26}}$, &\\$v_{23\overline{35}}$,$v_{33\overline{36}}$,$v_{\overline{1268}}$,$v_{\overline{1358}}$,$v_{\overline{1567}}$,$v_{\overline{2467}}$, &\\ $v_{\overline{3457}}$,$v_{1\overline{1246}}$,$v_{1\overline{1345}}$,$v_{2\overline{1256}}$,$v_{2\overline{2345}}$, &\\$v_{3\overline{1356}}$,$v_{3\overline{2346}}$,$v_{\overline{123456}}$ & $\omega_{14}=(a,b,c,t+6)$\\
	\hline 
	$v_{\overline{4567}}$,$v_{1\overline{1456}}$,$v_{2\overline{2456}}$,$v_{3\overline{3456}}$ & $\omega_{15}=(a,b,c+2,t+6)$\\
	\hline
\end{tabular}
\bigskip
\caption{Weights in degree six.}
\label{weightsDegree6}
\end{table}

\begin{lemma}\label{omega131415}
	The highest weight of $V$ is one among $\omega_{13}$, $\omega_{14}$, $\omega_{15}$.
\end{lemma}
\begin{proof}
The result follows from the same argument used in
the previous degrees, by considering the steps of
the weights and Table~\ref{weightsDegree6}.
\end{proof}
\begin{proposition}\label{omega14}
	$\omega_{14}$ is not the highest weight.
\end{proposition}
\begin{proof}
	Assume by contradiction that $\omega_{14}$ is the highest weight. Since $\omega_{15}$ is strictly higher than $\omega_{14}$, the four elements of weight $\omega_{15}$ vanish.
	If we assume $a\neq 0$ the equations given by $(x_ix_j\otimes z_\epsilon)w_{sing}(6)=0$ and $(d_{ij}\otimes z_\epsilon)w_{sing}(6)=0$ force all vectors of weight $\omega_{14}$ to vanish. By duality the same happens if $b\neq 0$ so we can assume $a=b=0$. In this case, repeated applications of
$x_ix_j\otimes z_\epsilon$ and
$d_{ij}\otimes z_\epsilon$, for $\epsilon=+,-$,
yield, after lengthy calculations, that all
vectors of weight $\omega_{14}$ vanish.
\end{proof}

Recall that we let $D_0=d_1^-d_2^-d_3^-$, $D_1=d_1^+d_2^-d_3^-+d_1^-d_2^+d_3^-+d_1^-d_2^-d_3^+$, $D_2=d_1^+d_2^+d_3^-+d_1^+d_2^-d_3^++d_1^-d_2^+d_3^+$ and $D_3=d_1^+d_2^+d_3^+$.
We also let, for $\epsilon\in \{+,-\}$, $\eta^\epsilon=\partial_1d_1^\epsilon+\partial_2 d_2^\epsilon+\partial_3 d_3^\epsilon$.

\begin{theorem}\label{ClassDeg6}
The following is a complete list, up to scalar factors, of singular vectors of degree $6$ and of the corresponding morphisms.
    
\begin{itemize}
		\item $w[6_{AC}]=D_3\big((\eta^--(c+1)d^-)\otimes z_+^{c-2}+(c-2)d^+\otimes z_+^{c-3}z_-\big)$ in $M_C(0,c-2)$, for all $c\geq 2$, inducing a morphism
		$\varphi[6_{AC}]:M_A(0,c)\rightarrow M_C(0,c-2)$;\vspace{2mm}
        \item $w[6_{BD}]=(\eta^++cd^+)\sum_{k=0}^3 D_k\otimes \de_+^k \de_-^{c+2-k}\in M_D(0,c+2)$ inducing a morphism \\$\varphi[6_{BD}]: M_B(0,c)\rightarrow M_D(0,c+2)$.
\end{itemize}
    
\end{theorem}
\begin{proof}
By Lemma \ref{omega131415}, Proposition \ref{omega14} and duality it is enough to show that if $\omega_{13}$ is the highest weight of $V$ then we have $(a,b,c,t)=(0,0,c,3c)$ for some $c\geq 2$ and $w_{sing}(6)$ is, up to a scalar factor, equal to $w[6_{AC}]$, and that $w[6_{BD}]$ is indeed a singular vector. The weights lower than $\omega_{13}$ among the $\omega_i$'s are $\omega_2, \omega_5, \omega_7$ (see Figure \ref{orderDegree6}) so we can assume that all the other weight spaces vanish. In particular $w_{sing}(6)$ is a sum of the 15 terms corresponding to the weights $\omega_2, \omega_5, \omega_7, \omega_{13}$.
 Applying $d_{12}\otimes z_-$ one immediately has $v_{1\overline{17}}=v_{2\overline{27}}=v_{3\overline{37}}=v_{2\overline{17}}=v_{3\overline{17}}=v_{\overline{1347}}=v_{\overline{1247}}=0$. We let $v=v_{1\overline{1234}}$. One also has $v_{2\overline{1235}}=v_{3\overline{1236}}=v_{\overline{2347}}=v_{\overline{1267}}=-v_{\overline{1357}}=v$. The equations $(a+b)v=0$ and $(-3c+t)v=0$ imply $a=b=0$ and $t=3c$ (otherwise all terms of weight $\omega_{13}$ vanish). We also have $v_{\overline{1238}}=(c+1)v$ and $v_{\overline{1237}}=Fv$ so that we can conclude that, up to a sign, $w_{sing}(6)=w[6_{AC}]$.
 
 The fact that $w[6_{BD}]$ is singular follows form a straightforward but tedious verification. 
\end{proof}

The following formulas will be useful in the sequel. Some of them follow immediately from definitions. Some others require more computations.
\begin{lemma}\label{relat} Let $D_{-1}=D_4=0$. Then, for
$k=0,1,2,3$ and $i=1,2,3$, we have:
\begin{itemize}
     \item $D_k d^\epsilon=-d^\epsilon D_k$ for $\epsilon=+,-$ and $k=0,1,2,3$;
     \item $d_i^+D_k+d_i^- D_{k+1}=0$ for $i=1,2,3$ and $k=0,1,2$;
     \item $D_k\eta^+ + \eta^+ D_k=kD_k d^++(k+1)D_{k+1}d^-$;
     \item $D_k \eta^- + \eta^- D_k=(4-k)D_{k-1}d^++(3-k)D_k d^-$;
     \item $[E,D_k]=(k+1)D_{k+1}$;
     \item $[F,D_k]=(4-k)D_{k-1}$;
     \item $[E,\eta^+]=0$, $[E,\eta^-]=\eta^+$;
     \item $[F,\eta^+]=\eta^-$, $[F,\eta^-]=0$;
     \item $[E,d^+]=0$, $[E,d^-]=d^+$;
     \item $[F,d^+]=d^-$, $[F,d^-]=0$.
 \end{itemize}
\end{lemma}
The following result expresses the two families of morphisms of degree 6 explicitly. Surprisingly, they can be obtained by the same formula with two different meanings.

% The inducing module in $M(0,0,{c},3 {c})$ is given by homogenoeus polynomials of degree ${c}$ in the variables $z_+$ and $z_-$ and, similarly, the inducing module in $M(0,0,{c}-2,3{c}+6)$ is given by homogeneous polynomials of degree ${c}-2$ in the variables $z_+$ and $z_-$. Therefore we denote $\partial_+$ and $\partial_-$ the derivatives with respect to these variables.

% We let 
% \[
% \eta_a=\partial_1d_1^a+\partial_2 d_2^a+\partial_3 d_3^a,
% \]
% for $a=+,-$.

\begin{theorem}\label{ClassDeg6Morphisms}The morphisms $\varphi[6_{AC}]$ and
$\varphi[6_{BD}]$ can both be expressed by the
following formula, with two different
interpretations of the operators involved:

\[
\varphi=d^-\sum_{k=0}^3 D_k\otimes \de_+^k\de_-^{3-k}z_{+}- d^+\sum_{k=0}^3 D_k\otimes \de_+^k\de_-^{3-k}z_{-}+ \eta^+ \sum_{k=0}^2 D_k \otimes \de_+^k \de_-^{2-k}.
\]

For the morphism $\varphi[6_{AC}]$, $z_{+}$ and $z_{-}$ act by multiplication, whereas $\de_+$ and $\de_-$ act by differentiation. For the morphism $\varphi[6_{BD}]$, the roles are reversed: $z_{+}$ and $z_{-}$ act by contraction, i.e.
\[
z_{\epsilon}\de_{\epsilon'}=-\delta_{\epsilon,\epsilon'},
\qquad \textrm{for }\epsilon,\epsilon'\in\{+,-\},
\]
whereas $\de_+$ and $\de_-$ act by multiplication.
\end{theorem}

%{\color{blue} In particular the corresponding singular vector of degree 6 is
%\begin{align*}
%\frac{1}{{c}({c}-1)}\varphi(z_+^{{c}})&=({c}-2)d_1^+d_2^+d_3^+d^+z_+^{{c}-3}z_-\\
%&-({c}+1)d_1^+d_2^+d_3^+d^-z_+^{{c}-2}\\
%&+(d_1^+d_2^+d_3^-+d_1^+d_2^-d_3^++d_1^-d_2^+d_3^+)\otimes z_+^{{c}-2}\\
%&-\eta_-d_1^+d_2^+d_3^+\otimes z_+^{{c}-2}
%\end{align*}}

 \color{black}
 \begin{proof}
 We first observe that $\varphi(\de_-^c)=w[6_{BD}]$, since $\eta^+D_3=0$, and  $\varphi(z_+^c)=c(c-1)w[6_{AC}]$, by Lemma \ref{relat}. 

 By Proposition~\ref{morphisms}, it is enough to
verify that $\varphi$ is a $\mathfrak g_0$-module
homomorphism. It is clear that $\varphi$ commutes
with the actions of $\mathfrak{sl}_3$, $Y$ and $H$.
Thus, it remains to verify that
\[
[E,\varphi]=[F,\varphi]=0.
\]

The identity $[E,\varphi]=0$ follows from
Lemma~\ref{relat}, using the actions of $E$ on
$D_k$, $\eta^+$, $d^+$ and $d^-$.
The identity $[F,\varphi]=0$ follows similarly,
using the additional relation
\[
\eta^-D_k+\eta^+D_{k-1}=0,
\]
which follows from the relations
$d_i^+D_k+d_i^-D_{k+1}=0$
in Lemma~\ref{relat}.

 \end{proof}
\section{Degree 7 and completion of the proof}

We have the following result for singular vectors of degree 7.
\begin{theorem}\label{ClassDeg7}
The following is a complete list, up to a scalar factor, of singular vectors of degree 7 and corresponding morphisms.
\begin{itemize} 
\item $w[7_{AD}]=D_3(\eta^--2d^-)\sum_{k=1}^3d_k^-\otimes\de_k$ in
$M_D(1,0)$ inducing a morphism\\ 
$\varphi[7_{AD}]:
M_A(0,1)\rightarrow M_D(1,0)$;\vspace{2mm}
\item $w'[7_{AD}]=\eta^+ d_1^+(D_1\otimes \de_++D_0\otimes \de_-)$  in $M_D(0,1)$ inducing a morphism\\
$\varphi'[7_{AD}]:
M_A(1,0)\rightarrow M_D(0,1)$.
\end{itemize}
    \end{theorem}
\begin{proof}We proceed as in the previous degrees and only
sketch the proof. The basic linear homogeneous system reduces the form of a generic singular vector to a sum of 29 summands. If we denote by $(a,b,c,t)$ the weight  of the singular vector, then  the coefficients of these 29 vectors have 4 possible weights: $\omega_1=(a-1,b,c-1,t+7)$, $\omega_2=(a-1,b,c+1,t+7)$, $\omega_3=(a,b+1,c-1,t+7)$ and $\omega_4=(a+1,b-1,c-1,t+7)$. By considering the possible weight steps, we see
that the highest weight must be either $\omega_2$
or $\omega_3$. These two weights are dual to each other. Assume that $\omega_2$ is the highest weight. In this case, all terms of weights $\omega_3$ and $\omega_4$ vanish,
since neither of these weights is dominated by
$\omega_2$. This reduces the singular vector to the sum of 19 terms.
If $a\neq 1$ or $b\neq 0$, applying the operators
$x_ix_j\otimes z_\epsilon$ shows that all
coefficients of weight $\omega_2$ vanish.
Thus, we may assume that $a=1$ and $b=0$. The singular vector at this stage is a sum of 16 terms in $M(0,0,c+1,t+7)$. Applying also elements $d_{ij}\otimes z_{\epsilon}$ we deduce that $c=0$ and $t=-4$, i.e. the singular vector has weight $(1,0,0,-4)$ and lies in $M(0,0,1,3)=M_D(0,1)=M_A(0,1)$. The same equations force the singular vector to be $w'[7_{AD}]$ in the statement. By duality, there is a unique other singular vector, up to scalar factors, and one can check that such vector is indeed $w[7_{AD}]$.
\end{proof}
\begin{theorem}\label{list}
    We have 
    \begin{align*}
    \varphi[7_{AD}]&=D_3(\eta^--2d^-)\sum_kd_k^-\otimes \de_k\de_+ + D_2(\eta^--2d^-)\sum_kd_k^-\otimes \de_k\de_-,\\
    \varphi'[7_{AD}]&=\eta^+\sum_k d_k^+ (D_1\otimes \de_++D_0 \otimes \de_-) \de_k.
    \end{align*}
\end{theorem}
\begin{proof}
The result follows by a direct verification using
Proposition~\ref{morphisms}. For $\varphi[7_{AD}]$,
the required identities follow from Lemma~\ref{relat}.

    The identity $[E,\varphi'[7_{AD}]]=0$ follows from
Lemma~\ref{relat}, together with the observation
that
\[
\eta^+d_k^+D_2=0
\qquad\text{for every }k=1,2,3.
\]
Similarly, the identity $[F,\varphi'[7_{AD}]]=0$
follows from Lemma~\ref{relat} and the relation
\[
\eta^-d_k^+D_1+\eta^+d_k^+D_0=0.
\]
\end{proof}
We have classified all singular vectors of degree $1,2,3,4,6,7$ and related morphisms between Verma modules. 
\begin{theorem}\label{FinalClass}Let $M(a,b,c,t)$ be a finite Verma module and let
$w\in M(a,b,c,t)$ be a nonzero homogeneous weight
vector of positive degree. Then $w$ is a singular
vector if and only if, up to a scalar factor,
$w$ is one of the following:

\medskip

\begin{enumerate}
    \item[]\hspace{-10mm}\emph {-- Degree 1:} $w[1_A]\in M_A(a-1,c-1)$, $w[1_B]\in M_B(a-1,c+1)$, $w[1_C]
\in M_C(b+1,c-1)$,  $w[1_D]\in M_D(b+1,c+1)$ as described in Theorem \ref{classDegree1};
\item[]\hspace{-10mm}\emph {-- Degree 2:} $w[2_{AB}]\in M_B(a-2,0)$, $w[2_{AD}]\in  M_D(1,2)$, $w[2_{CD}]\in M_D(b+2,0)$ as described in Theorem \ref{degree2};
\item[]\hspace{-10mm}\emph {-- Degree 3:} $w[3_{AC}]\in  M_C(0,c-3)$, $w[3_{BD}]\in M_D(0,c+3)$  as described in Theorem \ref{ClassDeg3};
\item[]\hspace{-10mm}\emph {-- Degree 4:}  $w[4_{AD}] \in M_D(1,0)$, $w'[4_{AD}]\in M_D(0,2)$ as described in Theorem \ref{ClassDeg4};
\item[]\hspace{-10mm}\emph {-- Degree 6:} $w[6_{AC}]\in M_C(0,c-2)$, $w[6_{BD}]\in M_D(0,c+2)$ as described in Theorem \ref{ClassDeg6};
\item[]\hspace{-10mm}\emph {-- Degree 7:} $w[7_{AD}]\in M_D(1,0)$ and $w'[7_{AD}]\in M_D(0,1)$ as described in Theorem \ref{ClassDeg7}.
\end{enumerate}
\end{theorem}
\begin{proof}By Corollary~\ref{bound}, the degree of a singular
vector is at most $14$. A computer verification
shows that the basic homogeneous linear system has
no nonzero solutions for
$d=5,8,9,10,11,12,13,14$.
The result then follows from
Theorems~\ref{classDegree1}, \ref{degree2},
\ref{ClassDeg3}, \ref{ClassDeg4},
\ref{ClassDeg6} and \ref{ClassDeg7}.
\end{proof}

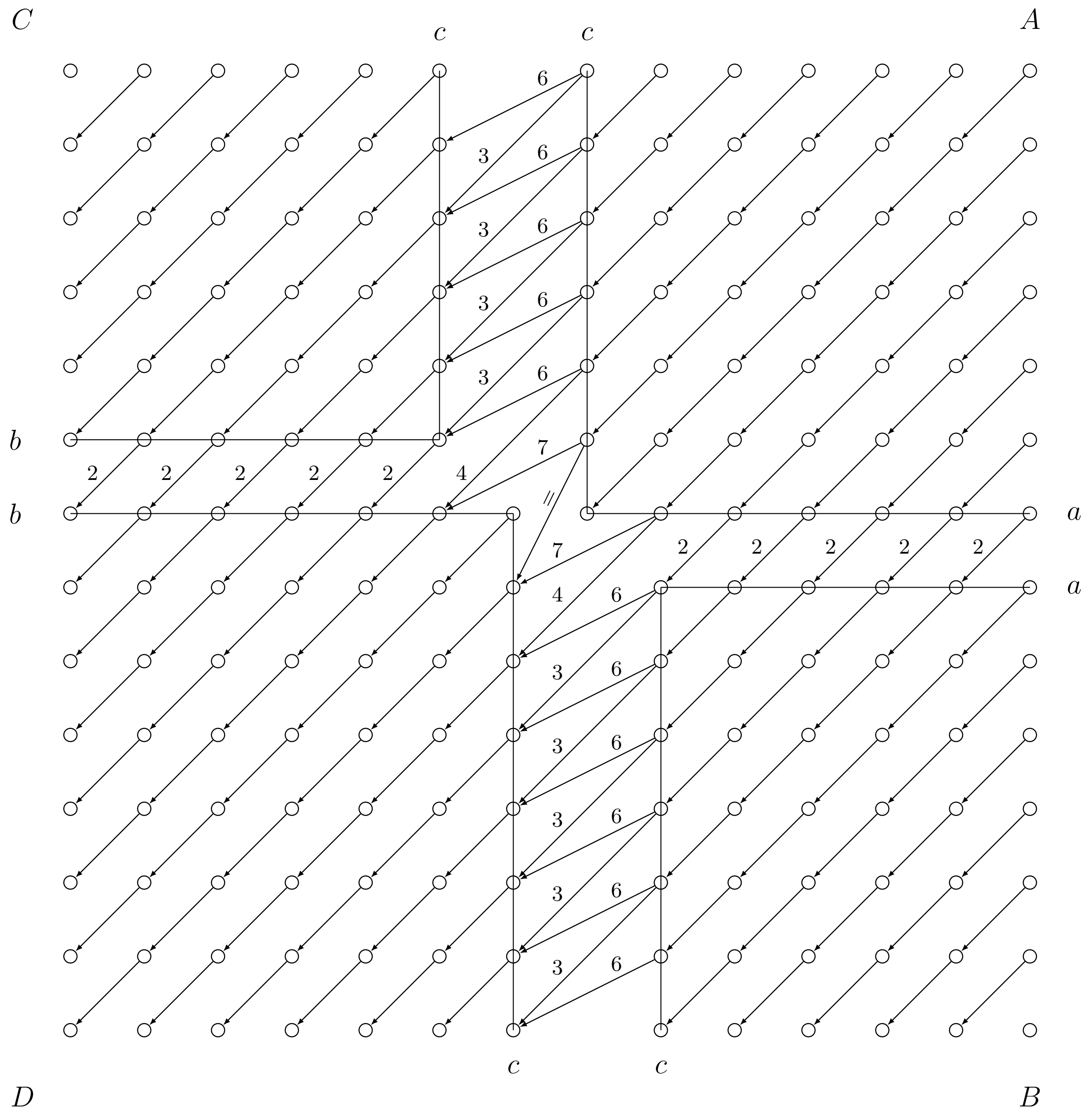
\begin{figure}
\begin{tikzpicture}[scale=1.5,
	x=0.72cm,
	y=0.72cm,
	>=Latex,
	dot/.style={
		circle,
		draw,
		fill=white,
		minimum size=5.5pt,
		inner sep=0pt,
		line width=.45pt
	},
	arr/.style={
		-{Latex[length=3pt,width=2.2pt]},
		line width=.45pt
	},
	darr/.style={
		-{Latex[length=3pt,width=2.2pt]},
		dash pattern=on 4pt off 4pt,
		line width=.55pt
	},
	barr/.style={
		-{Latex[length=4pt,width=3pt]},
		line width=1pt
	}
	]
	
	% =========================================================
	% QUADRANTE C
	% =========================================================
	
	\foreach \i in {0,...,5}{
		\foreach \j in {0,...,5}{
			\node[dot] (C\i\j) at (\i,\j) {};
		}
	}
	
	% frecce interne: da NE verso SW
	\foreach \i in {1,...,5}{
		\pgfmathtruncatemacro{\im}{\i-1}
		\foreach \j in {1,...,5}{
			\pgfmathtruncatemacro{\jm}{\j-1}
			\draw[arr] (C\i\j) -- (C\im\jm);
		}
	}
	
	\draw (C00.center) -- (C50.center);
	\draw (C50.center) -- (C55.center);

	% =========================================================
	% QUADRANTE A
	% =========================================================
	
	\foreach \i in {0,...,6}{
		\foreach \j in {0,...,6}{
			\node[dot] (A\i\j) at ({7+\i},\j-1) {};
		}
	}
	
	% frecce interne: da NE verso SW
	\foreach \i in {1,...,6}{
		\pgfmathtruncatemacro{\im}{\i-1}
		\foreach \j in {1,...,6}{
			\pgfmathtruncatemacro{\jm}{\j-1}
			\draw[arr] (A\i\j) -- (A\im\jm);
		}
	}
	
	\draw (A00.center) -- (A60.center);
	\draw (A00.center) -- (A06.center);

	% =========================================================
	% QUADRANTE D
	% =========================================================
	
	\foreach \i in {-1,...,5}{
		\foreach \j in {0,...,7}{
			\node[dot] (D\i\j) at (\i+1,-1-\j) {};
		}
	}
	
	% frecce interne
	\foreach \i in {0,...,5}{
		\pgfmathtruncatemacro{\im}{\i-1}
		\foreach \j in {0,...,6}{
			\pgfmathtruncatemacro{\jm}{\j+1}
			\draw[arr] (D\i\j) -- (D\im\jm);
		}
	}
	
	\draw (D50.center) -- (D57.center);
	\draw (D-10.center) -- (D50.center);

	% =========================================================
	% QUADRANTE B
	% =========================================================
	
	\foreach \i in {0,...,5}{
		\foreach \j in {0,...,6}{
			\node[dot] (B\i\j) at ({8+\i},-2-\j) {};
		}
	}
	
	% frecce interne
	\foreach \i in {1,...,5}{
		\pgfmathtruncatemacro{\im}{\i-1}
		\foreach \j in {0,...,5}{
			\pgfmathtruncatemacro{\jm}{\j+1}
			\draw[arr] (B\i\j) -- (B\im\jm);
		}
	}
	
	\draw (B00.center) -- (B50.center);
	\draw (B00.center) -- (B06.center);

	% =========================================================
	% FRECCE TRA QUADRANTI
	% =========================================================
	
	% A -> C
	
	\draw[arr] (A02) -- (C50);
	\draw[arr] (A03) -- (C50);
	\draw[arr] (A03) -- (C51);
	\draw[arr] (A04) -- (C51);
	\draw[arr] (A04) -- (C52);
	\draw[arr] (A05) -- (C52);
	\draw[arr] (A05) -- (C53);
	\draw[arr] (A06) -- (C53);
	\draw[arr] (A06) -- (C54);
	% B -> D
	\draw[arr] (B00) -- (D52);
	\draw[arr] (B00) -- (D53);
	\draw[arr] (B01) -- (D53);
	\draw[arr] (B01) -- (D54);
	\draw[arr] (B02) -- (D54);
	\draw[arr] (B02) -- (D55);
	\draw[arr] (B03) -- (D55);
	\draw[arr] (B03) -- (D56);
	\draw[arr] (B04) -- (D56);
	\draw[arr] (B04) -- (D57);
	\draw[arr] (B05) -- (D57);
	
	%C-->D
	\draw[arr] (C50) -- (D30);
	\draw[arr] (C40) -- (D20);
	\draw[arr] (C30) -- (D10);
	\draw[arr] (C20) -- (D00);
	\draw[arr] (C10) -- (D-10);
	
	%A-->B
	\draw[arr] (A20) -- (B00);
	\draw[arr] (A30) -- (B10);
	\draw[arr] (A40) -- (B20);
	\draw[arr] (A50) -- (B30);
	\draw[arr] (A60) -- (B40);

	% A -> D
        \draw[arr] (A10) -- (D51);
        \draw[arr] (A01) -- (D40);
	\draw[arr] (A02) -- (D40);
	\draw[arr] (A10) -- (D52);
	\draw[arr] (A01) -- (D51);
	
	% =========================================================
	% ETICHETTE
	% =========================================================
	
	\node at (-0.65,5.7) {$C$};
	\node at (13.0,5.7) {$A$};
	
	\node at (-0.65,-8.9) {$D$};
	\node at (13.0,-8.9) {$B$};
	
	\node at (5,5.5) {$c$};
	\node at (7,5.5) {$c$};
	
	\node at (-0.75,0) {$b$};
	\node at (-0.75,-1) {$b$};
	
	\node at (13.6,-1) {$a$};
	\node at (13.6,-2) {$a$};
	
	\node at (6,-8.5) {$c$};
	\node at (8,-8.5) {$c$};
	
	% etichette gradi
	\foreach \i in {1,...,5}{
		\node at (\i-0.7,-0.45) {\scriptsize{2}};
	}

	\foreach \i in {1,...,5}{
		\node at (\i+7.3,-1.45) {\scriptsize{2}};
	}
	
	\foreach \j in {1,...,4}{
		\node at (5.6,\j-0.15) {\scriptsize{3}};
	}
	
	\foreach \j in {1,...,5}{
		\node at (6.4,\j-0.1) {\scriptsize{6}};
	}

        \node at (6.4, -0.1) {\scriptsize{7}};
        \node at (6.6, -1.5) {\scriptsize{7}};
    
	\foreach \j in {1,...,5}{
		\node at (6.6,-\j-2.15) {\scriptsize{3}};
	}
	
	\foreach \j in {1,...,6}{
		\node at (7.4,-\j-1.1) {\scriptsize{6}};
	}
	
	\node at (5.3,-0.45) {\scriptsize{4}};
	\node at (6.6,-2.1) {\scriptsize{4}};
	\node[rotate=64] at (6.5,-0.8){\scriptsize{=}};
	
        \node at (2.6,-10) {$M_A(a,c)=M(a,0,c,-4a+3c)$;};
	\node at (9,-10) {$M_B(a,c)=M(a,0,c,-4a-3c-6)$;};
    
%	\node at (5,-10) {$M_A(a,c)=M(a,0,c,-4a+3c)\hspace{8mm}$};
%	\node at (5,-11) {$M_B(a,c)=M(a,0,c,-4a-3c-6\hspace{1.5mm})$};
	\node at (2.8,-11) {$M_C(b,c)=M(0,b,c,4b+3c+12)$;};
	\node at (8.8,-11) {$M_D(b,c)=M(0,b,c,4b-3c+6)$.};
	
\end{tikzpicture}
\caption{\label{FinalFigure} All morphisms between finite Verma modules}
\end{figure}

The classification in Theorem~\ref{FinalClass}
is summarized in Figure~\ref{FinalFigure}.
Every morphism is represented by an arrow labelled
with its degree, with labels equal to $1$ omitted.
The only morphism not explicitly represented is
$\varphi[2_{AD}]$, which is the composition of two
degree-$1$ morphisms, as described in
Theorem~\ref{degree2}.

\bigskip

\begin{corollary}\label{V.}
All sequences in Figure \ref{FinalFigure} are complexes except for the sequence passing through $M_A(0,1)=M_D(0,1)$. The composition of morphisms $M_A(1,2)\rightarrow M_A(0,1)$ and $M_D(0,1)\rightarrow M_D(2,1)$ is not 0, but if we replace these two morphisms of degree 1 by their composition, the sequence become a complex.
\end{corollary}

\end{document}